\documentclass[11pt,reqno]{amsart}

\usepackage[x11names]{xcolor}
\usepackage[linkcolor=Blue2,citecolor=black,colorlinks]{hyperref}
\usepackage{color, bm, amscd, tikz-cd}

\usepackage{amsmath,amsthm,amssymb, mathrsfs, bbm}

\allowdisplaybreaks

\theoremstyle{plain}
\newtheorem{theorem}{Theorem}[section]

\newtheorem{lemma}[theorem]{Lemma}
\newtheorem{proposition}[theorem]{Proposition}

\theoremstyle{definition}
\newtheorem{definition}[theorem]{Definition}

\newtheorem{remark}[theorem]{Remark}

\numberwithin{equation}{section}

\newcommand{\R}{\mathbb{R}}
\newcommand{\C}{\mathbb{C}}

\newcommand{\D}{\mathbb{D}}
\newcommand{\T}{\mathbb{T}}

\newcommand{\bydef}{\stackrel{\rm def}{=}}
\usepackage{mismath}

\begin{document}
\title[Superresolution]{The multivariate Herglotz-Nevanlinna class: superresolution}

\author[M. Bhowmik]{Mainak Bhowmik }
\address{Indian Institute of Science, Bangalore, India }
\email{\tt mainak.bhowmik943@gmail.com, mainakb@iisc.ac.in}

\author[M. Putinar]{Mihai Putinar}
\address{University of California at Santa Barbara, CA} 
\email{\tt mputinar@math.ucsb.edu}


\keywords{superresolution, Herglotz-Nevanlinna class, positive pluriharmonic function, 
Carath\'eodory-Fej\'er interpolation, Nevanlinna-Pick interpolation, rational inner function, polydisk}

\subjclass[2020]{32E30, 31B20, 32A17, 30E05, 44A60} 
\dedicatory {Bent Fuglede, in memoriam}

\begin{abstract} Bounded holomorphic interpolation problems associated to finitely many data have, in general, distinct solutions. Uniqueness arises only in some convex extreme configurations.
Rational inner functions in a polydisk are the best understood examples in this sense. We analyze the continuity of global solutions as functions of the finite interpolation data in neighbourhoods of elements distinguished by this uniqueness property. Our study covers rational inner or Cayley rational inner functions in
	the polydisk and automorphisms of the Euclidean ball. The proof of the main superresolution result is derived from optimization theory techniques and volume estimates of sublevel sets of real polynomials, both emerging from Markov's multivariable moment problem.

\end{abstract}

\maketitle

\section{Introduction} This is the first part of a work devoted to the
Herglotz-Nevanlinna class of functions depending on several complex variables.
The concept of superresolution is closely related to that of spectrum. The former covers several meanings, all related to specific inverse problems. Treating 
the spectrum as a "ghost" of an inaccessible entity, a widely accepted 
interpretation of the superresolution phenomenon is the effective approximation of
the entire spectrum from a part of it in the proximity of a finitely determined locus. Notable recent advances in signal/image processing validate this point of view \cite{Candes-Fernandez, Donoho, Mallat}.

A more precise manifestation of superresolution was proposed in the setting of (truncated) moment problems, specifically the Markov $L$-problem \cite{Krein-Nudelman, Putinar-JAT}. An entropy optimization problem naturally entering into these studies
 \cite{Budisic-Putinar, Gamboa, Lewis, Putinar}. 
The present note is a continuation of these investigations in the framework of bounded holomorphic interpolation of functions depending on
several complex variables. Along these lines, the gap between a single complex variable and two or more complex variables is sizable. We elaborate below some 
observations regarding the multivariable setting, which in their turn 
raise challenging open problems.

The class of holomorphic functions in the unit disk $\D$ of the complex plane $\C$ with non-negative real part is parametrized via the classical Riesz-Herglotz integral transform by all positive Borel measures supported on the unit circle. A Cayley transform establishes a bijective correspondence of this class to the set of all holomorphic maps from the unit disk to itself. For the later class, a continued fraction type expansion combined with Schwarz Lemma led Schur to an algorithm bearing his name. Schur's algorithm provides a parametrization of holomorphic maps from the disk to itself which involves simple algebraic operations on the Taylor coefficients.
The algorithm finishes in finitely many steps if and only if the original is a rational inner function (i.e., a rational function with poles off $\D$ and having modulus equal to $1$ on the boundary). Transferring the above framework to holomorphic functions mapping the upper-half plane to itself brings into the picture the resolvents of (possibly unbounded) self-adjoint operators acting on a Hilbert space. Then naturally, it includes the power moment problem 
on the real line and the theory of orthogonal polynomials,
with all ramifications. Nevanlinna's continued fraction expansion for these functions paved the road of a century of surprising discoveries, not least touching the applied sciences realm. The mathematical construct emerging from this interaction of function theory of complex variables, Hilbert space operators and constructive approximation has reached maturity a few decades ago and continues to fascinate generations of researchers \cite{Agler-McCarthy-Young, Ahiezer-Krein, Krein-Nudelman}.

In spite of sustained theoretical efforts and targeted applications demands, the analogous structure of holomorphic functions defined on classical domains of $\C^d$ ($d>1$) with non-negative real part remains fragmented.
While a corresponding Riesz-Herglotz integral formula is true in any dimension, the positive measures associated to these functions are restricted by constraints naturally imposed by the boundary values of positive pluriharmonic functions. The operator theory counterpart is also confined to the smaller class of Herglotz-Agler functions. The next section contains precise definitions and bibliographical references.

{One of the better-understood situations is offered by the polydisk as a geometric support. This is due to the fact that rational inner functions (refer to Definition \ref{def:RIF}) abound: they interpolate finitely many Taylor coefficients (in 2D) and approximate, with respect to the Fr\'echet topology of locally uniform convergence, holomorphic functions bounded by one. Rational inner functions in the polydisk are determined by a finite section of their 
Taylor series at the origin or values on a prescribed finite set; cf. Pfister \cite[Satz 3]{Pfister}. Cayley rational inner functions (see Definition \ref{def:RCIF}) also exhibit a similar finite determinateness property.
To be more precise, let $\T$ denote the torus, seen as the boundary of the disk $\D$ and let $\mathbb{C}_{+}$ denote the right half plane in $\mathbb{C}$. A rational
function with poles off the polydisk $ R: \D^d \longrightarrow \mathbb{C}_{+}$ which maps $\T^d$ to the imaginary axis almost everywhere is determined by a finite section $T_{\boldsymbol{n}}(R) = \sum_{\boldsymbol{\alpha} \leq \boldsymbol{n}} c_{\boldsymbol{\alpha}}(R) \boldsymbol{z}^{\boldsymbol{\alpha}}$ of its Taylor expansion at zero:
$$ R(\boldsymbol{z}) = T_{\boldsymbol{n}}(R) +  \sum_{\boldsymbol{\alpha} \nleq \boldsymbol{n}} c_{\boldsymbol{\alpha}}(R) \boldsymbol{z}^{\boldsymbol{\alpha}}, \ \ \boldsymbol{z} \in \D^d.$$
}

{It is on the polydisk and the non-generic collection of Cayley rational inner functions that we illustrate a superresolution phenomenon.
We prove that, for any holomorphic function $f: \D^d \longrightarrow \mathbb{C}_{+}$ belonging to a fixed neighbourhood of $R$, the distance from $f$ to $R$ in the metric topology of ${\mathcal O}(\D^d)$, the collection of all holomorphic functions on $\D^d$, is bounded by a fractional power of the Euclidean distance between $T_{\boldsymbol{n}}(R)$ and $T_{\boldsymbol{n}}(f)$. Notations and precise bounds are contained in the subsequent sections.
We reduce the proof to the known setting of power moments of a shade function (that is a measurable function with values in $[0,1]$) defined on the torus $\T^d$. It was recently demonstrated in \cite{Putinar} that "black and white pictures" (that is characteristic functions) delimited by single polynomial equations are not only finitely determined by their moments, but they possess a superresolution property. That is, any shade close enough to them with respect to the finite
determining set of moments is globally close (in a precise sense) to the original. Markov's logarithmic transform of a Herglotz-Nevanlinna function in the polydisk makes this reduction possible \cite{Budisic-Putinar}.
}

{
The contents is the following.
Section \ref{Prelim} is devoted to some preliminaries on function theory in the polydisk.
The main superresolution result is contained in Section \ref{Main results}, covering small perturbations of data associated with a Cayley rational inner function in the polydisk. Then a simpler superresolution theorem is derived for small perturbations of the affine Taylor expansion of holomorphic automorphisms of the Euclidean unit ball. 
At the end, we provide some examples exhibiting non-uniqueness of solutions to finite bounded interpolation data and formulate an intriguing open question. }
\bigskip

{\bf Acknowledgements.} Both authors were able to visit each other's institutions and start a collaboration thanks to the generous support provided by the SPARC (Project No. 1698) travel grant, Govt. of India. The first author was also supported by the Prime Minister’s Research Fellowship PMRF-21-1274.03. The second author was partially supported by a Simons Collaboration Grant for Mathematicians. We thank the referee for high-quality comments.

\section{Herglotz-Nevanlinna class in the polydisk} \label{Prelim}

The present section aims to recall some terminology and collect a series of known ingredients that will appear in the body of the note.

\begin{definition}
A holomorphic function $f: \D^d \rightarrow \mathbb{C}$ is said to be a 
\begin{enumerate}
\item[(i)] {\em Schur function}, provided $|f(\boldsymbol{z})| \leq 1$ for every $\boldsymbol{z} \in \D^d$;

\item[(ii)] {\em Herglotz-Nevanlinna function}, provided $\Re f(\boldsymbol{z}) \geq 0$ for every $\boldsymbol{z} \in \D^d$.
\end{enumerate}
 The set of all Schur functions on $\D^d$ is denoted by $\mathcal{S}(\D^d)$ and the set of all Herglotz-Nevanlinna functions is $\mathcal{H}(\D^d)$. 
\end{definition}

We warn the reader on a lack of consensus concerning terminology arising from the passage to the multivariable setting. For some authors, Nevanlinna functions are
holomorphic functions with non-negative imaginary part defined on the poly-upper half plane. Others call them Herglotz-Nevanlinna functions or even Herglotz functions for short. Also, the term
Herglotz function is sometimes attached to holomorphic functions with non-negative real part defined on a given domain. Throughout this article we will make precise what class of functions is involved in the
statements and proofs.

Schur functions in one dimension are well understood from a complex analytic as well as an operator theoretic perspective. Among various characterization of a Schur function none is more fascinating than the so called {\em Schur's algorithm} which detects contractivity of a holomorphic function in $\D$ through countably many complex numbers known as the {\em Schur parameters}. Following \cite{Khrushchev}, we describe this algorithm as a prelude to illustrate the superresolution phenomenon by leveraging the information encoded within the Schur parameters.

\subsection{Schur's algorithm}
Let $f$ be a fixed Schur function on $\D$. Let $f_0 \bydef f$ and $\gamma_0(f) \bydef f_0(0)$. For $n \in \mathbb{N}_0$, we define inductively
\begin{align} \label{Eq:Schur parameter}
\gamma_n \bydef f_n(0), \, \, f_{n+1}(z) \bydef \frac{f_n(z)- \gamma_n}{z(1-\overline{\gamma_n} f_n(z))}.
\end{align}
 Note that Schwarz Lemma guarantees that each $f_n$ defines a Schur function on $\D$. Moreover, the above relations can be rewritten as 
 \begin{align}\label{Eq:Schur inversion}
 f_n(z)= \gamma_n +  \frac{(1-|\gamma_n|^2)z}{\overline{\gamma_n} z + 1/f_{n+1}(z)}.
 \end{align}
Iterating the above expression, we obtain the so-called Wall continued fraction:
\begin{align} \label{Eq:contd frac}
f(z)= \gamma_0 + \frac{(1-|\gamma_0|^2)z}{\overline{\gamma_0} z}_ + \frac{1}{\gamma_1}_+ \frac{(1-|\gamma_1|^2)z}{\overline{\gamma_1} z}_ {+\dots}.
\end{align}
The complex numbers $\gamma_n$ are called the Schur parameters of $f$. The algorithm continues (possibly indefinitely) until a first parameter $\gamma_k$ attains modulus one: $|\gamma_k|=1$. If $|\gamma_k|=1$ then by the maximum modulus principle $f_k \equiv \gamma_k$ and hence the algorithm breaks down. Accordingly, we put $\gamma_n =1$ for $n>k$. Moreover, in this case the function $f$ must be a Blaschke product of degree $k$, that is a finite product of M\"obius transforms.

A polynomial (in one variable) $p$ of degree $n$ gives rise to its {\em reflection polynomial} $p^*(z) \bydef z^n \overline{p(\frac{1}{\bar{z}})}$. Later on we shall define this notion for multivariable polynomial as well. Clearly, the degree of $p^*$ is at most $n$. Now, we inductively define the Wall polynomials related to the continued fraction \eqref{Eq:contd frac} as follows:
\begin{align*}
A_0 \equiv \gamma_0, A_0^*\equiv \bar{\gamma}_0, B_0= B_0^* \equiv 1; \\
A_{n+1}(z)= A_n(z)+ \gamma_{n+1}zB_n^*(z), & \, \, B_{n+1}(z)= B_n(z)+ \gamma_{n+1}zA_n^*(z), \\
A_{n+1}^*(z)=zA_n^*(z)+ \bar{\gamma}_{n+1} B_n(z), & \, \, B_{n+1}^*(z)=zB_n^*(z)+ \bar{\gamma}_{n+1} A_n(z).
\end{align*}
It is easy to see that $A_n$ and $B_n$ are both polynomials of degree $n$. The following proposition describes various properties of these polynomials.
\begin{proposition} \label{Prop:Wall poly} 
	
\begin{enumerate}
\item $B_n^* B_n - A_n^* A_n = z^n \prod_{j=1}^n (1-|\gamma_j|^2) \bydef  \omega_n z^n$.

\item The polynomials $B_n$ do not vanish in $\overline{\D}$.

\item The functions $A_n/B_n$ and $A_n^*/B_n$ are in $\mathcal{S}(\D)$. Moreover, for each $n$ with $|\gamma_n|<1$,
\begin{align*}
\sup_{z\in \overline{\D}} | A_n(z)/B_n(z)| < 1, \, \, \text{ and } \, \sup_{z\in \overline{\D} } | A_n^*(z)/B_n(z)| < 1.
\end{align*}  
\item For each $n$, we can express $f$ as 
$$f(z)= \frac{A_n + zB_n^* f_{n+1}}{B_n + zA_n^* f_{n+1}}.$$

\item The Schur function $f$ is a Blaschke product of degree $k$ if and only if $|\gamma_k|=1$. In that case, $f=\frac{A_k}{B_k}$.
\end{enumerate}
\end{proposition}

Part ($4$) of the above proposition implies that any Schur function $g$ on $\D$ whose first $n+1$ Schur parameters $\{\gamma_0, \dots, \gamma_n\}$ are fixed and $|\gamma_n|<1$, is of the form:
\begin{align*}
g(z)= \frac{A_n(z) + zB_n^*(z) h(z)}{B_n(z) + zA_n^*(z) h(z)} \, \, \text{ for some } h\in \mathcal{S}(\D).
\end{align*}
Moreover, if $\tilde{g}$ is a Schur function whose first $n+1$ Schur parameters are the same as that of $g$ then $g$ and $\tilde{g}$ have the same Taylor polynomial of order $n$ about the origin. The following remarkable characterization of functions in $\mathcal{S}(\D)$ is due to Schur.

\begin{theorem}[Schur \cite{Schur}]
The necessary and sufficient condition for a holomorphic function $f$ on $\D$ to be in $\mathcal{S}(\D)$ is that all  Schur parameters $\gamma_n \in \D$. Moreover, for given $\{\gamma_n \in \D: n\in \mathbb{N}_0\}$, we can reconstruct the Schur function via the continued fraction \eqref{Eq:contd frac} whose even convergents (i.e., convergents of order $2n$), $\frac{A_n}{B_n}$ converge uniformly on compact subsets of $\D$ to $f$.
\end{theorem}
An analogue of the Schur algorithm is unknown in higher dimensions. In an inspired early attempt to extend Schur's ideas to the multivariable setting,  Pfister \cite[Satz 14]{Pfister} has characterized the contractivity of a holomorphic function in $\D^d$ via the contractivity of infinitely many block Toeplitz matrices (of finite order) whose entries depend on the Taylor coefficients of the function under consideration.

\noindent \textbf{Multi-index notations:} The following notations are adopted for the rest of this note.
\begin{itemize}
\item[•] $\boldsymbol{\alpha} = (\alpha_1, \dots, \alpha_d) \in \mathbb{N}_0^d \text{ or } \mathbb{Z}^d$;
\item[•] $\boldsymbol{z}=(z_1, \dots, z_d) ,\ \ \boldsymbol{z}^{\boldsymbol{\alpha}} = \prod_{j=1}^d z_j^{\alpha_j}, \ \ \frac{1}{\overline{\boldsymbol{z}}} =(\frac{1}{\bar{z}_1}, \dots, \frac{1}{\bar{z}_d}) $;
\item[•] $\Gamma_{\boldsymbol{n}} = \{\boldsymbol{\alpha} \in \mathbb{N}_0^d: 0 \leq \alpha_j \leq n_j \ \text{ for } 1\leq j \leq d \}$ for $\boldsymbol{n}=(n_1, \dots, n_d)$ in $\mathbb{N}_0^d$;
\end{itemize}

\subsection{Representation of Herglotz-Nevanlinna functions}
For a Herglotz-Nevanlinna function $\varphi$ on $\D$, the function $f= (\varphi -1)/((\varphi +1))$ belongs to $\mathcal{S}(\D)$. This way, a bijective correspondence between $\mathcal{H}(\D)$ and $\mathcal{S}(\D)$ is established. The famous Riesz-Herglotz representation associates $\varphi$ with a unique positive regular Borel measure $\mu$ on the unit circle $\mathbb{T}$ via the integral formula:
\begin{align*}
\varphi(z)= i \Im \varphi(0) + \int_{\mathbb{T}} \frac{\xi + z}{\xi -z} d\mu(\xi) \, \, \text{ for } z\in \D.
\end{align*}
This integral representation also connects the moment problem with Schur's algorithm. In the polydisk $\D^d$, a version of the above integral representation was obtained by Kor\'anyi and Puk\'anszky \cite{KP}. It can be described as the following: for $\varphi$ in $\mathcal{H}(\D^d)$, there exists a unique positive regular Borel measure $\nu$ on the distinguished boundary $\mathbb{T}^d$ of $\D^d$ such that
\begin{align*} 
			\varphi(\boldsymbol{z})= i\Im \varphi(0)+ \int_{\mathbb{T}^d} H(\boldsymbol{z}, \boldsymbol{\xi}) d\nu(\boldsymbol{\xi}),
\end{align*}
where $$ H(\boldsymbol{z}, \boldsymbol{\xi}) = \frac{2}{\prod_{j=1}^{d} (1 - z_j \bar{\xi}_j)}-1, $$ and $\nu$ satisfies the moment conditions,
	\begin{equation} \label{Eq:Pluri-measure} 
\widehat{\nu}(n_1, \dots, n_d)= \int_{\mathbb T^d} \bar{\xi}_1^{n_1} \cdots \bar{\xi}_d^{n_d} d\nu(\boldsymbol{\xi}) = 0,
\end{equation} 
unless $n_j \geq 0$ for all $j=1, \cdots, d$ or  $n_j \leq 0$ for all $j=1, \cdots, d$.

Thus the pluriharmonic function $\Re \varphi(\boldsymbol{z})$ can be written as 
\begin{align}\label{Eq:KP-rep}
\Re \varphi(\boldsymbol{z}) = \int_{\mathbb{T}^d} \mathcal{P}_{\D^d}(\boldsymbol{z}, \boldsymbol{\xi}) d\nu(\boldsymbol{\xi})
\end{align}
where $\mathcal{P}_{\D^d}$ is the Poisson-Szeg\"o kernel of $\D^d$:
$$
\mathcal{P}_{\D^d}(\boldsymbol{z}, \boldsymbol{\xi}) = \prod_{j=1}^{d} \frac{ (1-|z_j|^2)}{ |1 - z_j \bar{\xi}_j|^2}.
$$
The positive regular Borel measures satisfying the above moment conditions \eqref{Eq:Pluri-measure} are often called {\em pluriharmonic measures} as their Poisson-Szeg\"o integral produces pluriharmonic functions (that are locally real part of holomorphic functions). 
 
As a matter of fact, if $\varphi$ extends continuously on $\overline{\D^d}$ then the representing measure $\nu$ is absolutely continuous with respect to the normalized Haar measure, $\Theta_d$ (sometimes denoted by $\Theta$) on $\mathbb{T}^d$ and $\Re \varphi$ is the density of $\nu$. For a general Herglotz-Nevanlinna function $\varphi$, the radial limit 
$$
\lim_{r \to 1^-} \Re \varphi (r \boldsymbol{\xi})
$$
exists for $\Theta_d$-a.e. $\boldsymbol{\xi}$ on $\mathbb{T}^d$ and it is the density of the absolute continuous part of $\nu$ with respect to $\Theta_d$.

Recently, in \cite{Bhowmik-BLMS} an integral representation for Herglotz-Nevanlinna functions in the style of Kor\'anyi and Puk\'anszky has been demonstrated. Furthermore, the interaction between the Herglotz-Nevanlinna representation and the Carath\'eodory approximation in finitely connected planar domains has been explored. See also \cite{BK-PAMS2} for an operator theoretic approach to Herglotz-Nevanlinna representation in a broader context.

\subsection{Carath\'eodory-Fej\'er interpolation on $\D^d$}
Given a set of scalars 

\noindent $\{c_{\boldsymbol{\alpha}}: \boldsymbol{\alpha} \in \Gamma_{\boldsymbol{n}} \}$ (called as the interpolation data), the {\em Carath\'eodory-Fej\'er interpolation} in the Schur class on $\D^d$ asks for a function $f$ in $\mathcal{S}(\D^d)$ with the Taylor series around the origin
\begin{align*}
f(\boldsymbol{z}) = \sum_{\boldsymbol{\alpha} \in \mathbb{N}_0^d} a_{\boldsymbol{\alpha}} \boldsymbol{z}^{\boldsymbol{\alpha}} \, \text{ for } \boldsymbol{z} \in \D^d 
\end{align*}
such that $a_{\boldsymbol{\alpha}} = c_{\boldsymbol{\alpha}}$ for $\boldsymbol{\alpha} \in \Gamma_{\boldsymbol{n}}$.

The Carath\'eodory-Fej\'er interpolation in the Herglotz-Nevanlinna setting can be formulated as the following. Given the interpolation data 

\noindent $\{ c_{\boldsymbol{\alpha}} : \boldsymbol{\alpha} \in \Gamma_{\boldsymbol{n}} \}$, we seek for a function $\varphi$ in the Herglotz-Nevanlinna class $\mathcal{H}(\D^d)$ such that the Taylor coefficients of $\varphi$ at the origin indexed by $\boldsymbol{\alpha}$ in $\Gamma_{\boldsymbol{n}} $ are $c_{\boldsymbol{\alpha}}$.

In one dimension, contractivity of a Toeplitz matrix formed by the given data determines the solvability of the interpolation problem in the Schur class, while the positivity of a Hankel matrix with entries from the given data ensures solvability in the Herglotz-Nevanlinna class. See \cite{FF} for a detailed treatment. In higher dimension ($d>1$), the solvability of the interpolation problem exists, but it is rather technical \cite{Dautov-Khud, Eschmeier-Patton-Putinar, Woerdeman}.

\subsection{Finitely determined functions} \label{subsection:Finite-determinateness} {The key protagonists of this note are the holomorphic maps uniquely determined by a finite segment of their Taylor series. This subsection discusses functions that exhibit such a finite determinateness property}.

The Finite Blaschke products are rational functions with poles off $\overline{\D}$ and unimodular on the boundary $\mathbb{T}$ of $\D$. This leads to the notion of rational inner function. 
\begin{definition} \label{def:RIF}
A bounded holomorphic function $f$ on $\D^d$ is said to be rational inner provided $f$ is rational with poles off $\D^d$ and is unimodular $\Theta_d-$a.e. on $\mathbb{T}^d$.
\end{definition}
In other terms, finite Blaschke products are rational inner functions on $\D$, and a simple exercise in complex analysis shows that these are all rational inner functions in $\D$. The structure of rational inner function on $\D^d$ was obtained in \cite[Satz 2]{Pfister} and it is described below. For a polynomial $p$ in $\mathbb{C}[\boldsymbol{z}]$ of multi-degree $\boldsymbol{n}=(n_1, \dots, n_d)$, the reflection of $p$ is defined to be the polynomial
$$
p^*(\boldsymbol{z}) = \boldsymbol{z}^{\boldsymbol{n}} \overline{p \left(\frac{1}{\overline{\boldsymbol{z}}} \right)}.
$$

The multi-degree $\boldsymbol{n}=(n_1, \dots, n_d)$ of a polynomial $p$ is defined by $n_j = \deg_{z_j} p, \ \ 1 \leq j \leq d.$ Note that, the coefficient of $\boldsymbol{z}^{\boldsymbol{n}}$ can be null.

\begin{theorem}[Pfister \cite{Pfister}] \label{Th:Pfister-RIF-structure}
If $f$ is a rational inner on $\D^d$, then there exist a polynomial $p$ in $\mathbb{C}[\boldsymbol{z}]$ whose zeros are outside $\D^d$ and $\boldsymbol{m} \in \mathbb{N}_0^d$ such that 
$$ f(\boldsymbol{z}) = \boldsymbol{z}^{\boldsymbol{m}} \frac{p^*(\boldsymbol{z})}{p(\boldsymbol{z})}. $$
Conversely, given any polynomial $p$ in $\mathbb{C}[\boldsymbol{z}]$ having zeros outside $\D^d$, the holomorphic function $\frac{p^*}{p}$ is a rational inner function on $\D^d$.
\end{theorem}

For an authoritative account of rational inner functions in the polydisk we refer to the recent survey article by Knese \cite{Knese}. 
Due to the central role in function theory, the structure of rational inner functions has also been explored on various other domains such as the general bounded symmetric domains \cite[Theorem~ 3.3]{Koranyi-Vagi} or certain quotient domains related to the polydisk \cite[Theorem~ 4.5]{Bhowmik-Kumar-JFA}. 
For a long time, the existence of inner functions on the unit ball in $\C^d$ has circulated as an open problem \cite{Aleksandrov}.

According to Schur's algorithm, a Blaschke product of degree $n$ in $\D$ is determined by its first $n+1$ Schur parameters, which in their turn, solely depend on the first $n+1$ Taylor coefficients of the Blaschke product. Although in higher dimension the Schur algorithm is unavailable, Pfister \cite[Satz 3]{Pfister} established that {\em any rational inner function $\frac{p^*}{p}$ on $\D^d$ is determined by its Taylor polynomial of multi-degree $\boldsymbol{m}=(m_1, \dots, m_d)$ where $\boldsymbol{m}$ is the multi-degree of the polynomial $p$}.

In the Herglotz-Nevanlinna setting the {\em Cayley rational inner functions} are analogous to the rational inner functions in the Schur class. 
\begin{definition} \label{def:RCIF}
A holomorphic function $\varphi$ with non-negative real part on $\D^d$ is said to be Cayley rational inner if $\varphi$ is rational with poles off $\D^d$ and $\Re \varphi$ is zero $\Theta_d-$ a.e. on $\mathbb{T}^d$.
\end{definition}
It is easy to check that a Cayley rational inner function $\varphi$ is just the Cayley transform $(z \mapsto (1+z)/(1-z))$ of the rational inner function, $f = (\varphi -1)/(\varphi +1)$. Therefore, {\em $\varphi$ is also determined by its Taylor polynomial of multi-degree $\boldsymbol{m}=(m_1, \dots, m_d)$, the multi-degree of $\varphi$, as so is true for $f$}.

The Kor\'anyi-Puk\'anszky integral representation \eqref{Eq:KP-rep} for $\Re \varphi$ combined with the power series expansion of $\mathcal{P}_{\D^d}$ can be rewritten as 
$$
\Re \varphi(\boldsymbol{z}) = \sum_{\boldsymbol{\alpha}\in \mathbb{Z}^d} \widehat{\nu}(\boldsymbol{\alpha}) r_1^{|\alpha_1|} \dots r_d^{|\alpha_d|} e^{\iota \boldsymbol{\alpha}\cdot \boldsymbol{\theta}},
$$ 
where $\boldsymbol{z} = (r_1 e^{\iota \theta_1}, \dots, r_d e^{\iota \theta_d})$, $\boldsymbol{\alpha}\cdot \boldsymbol{\theta} = \sum_{j=1}^d \alpha_j \theta_j$ and $\nu$ satisfies the moment conditions as described in \eqref{Eq:Pluri-measure}. Thus the finite determinateness of $\varphi$ implies that the representing measure $\nu$ is determined by the finitely many moments $\{\widehat{\nu}(\boldsymbol{\alpha}): 0 \leq \alpha_j \leq m_j \text{ for } j=1, \dots, d \}$. This is one of the crucial observations for our superresolution result.

\subsection{The exponential transform} \label{subsection:Exp-transform}
Every Herglotz-Nevanlinna function $\varphi$ on $\D^d$ admits an holomorphic logarithm, $\log \varphi$ whose imaginary part (the phase) lies in the interval $[ -\pi/2, \pi/2 ]$. Then $\psi = \frac{1}{2} + \frac{\iota}{\pi} \log \varphi $ is a Herglotz-Nevanlinna function with 
$$
\Re \psi(\boldsymbol{z}) \in [0, 1] \ \text{ for } \boldsymbol{z}\in \D^d.
$$
The radial limit of the function $\Re \psi$ produces an integrable, positive and bounded function $g$ on $\mathbb{T}^d$.
Thus the measure $d\nu$ in the Kor\'anyi-Puk\'anszky representation \eqref{Eq:KP-rep} for $\Re \psi$ is $g d\Theta$. Hence,
\begin{align} \label{Eq:Exp}
\varphi(\boldsymbol{z})= \iota \exp (\pi \Im \psi(0))  \exp \left(-\iota \pi  \int_{\mathbb{T}^d} H(\boldsymbol{z}, \boldsymbol{\xi}) g(\boldsymbol{\xi})d\Theta(\boldsymbol{\xi}) \right)
\end{align}
and 
\begin{align*}
\Re \psi(\boldsymbol{z}) = \int_{\mathbb{T}^d} \mathcal{P}_{\D^d}(\boldsymbol{z}, \boldsymbol{\xi}) g(\boldsymbol{\xi})d\Theta(\boldsymbol{\xi}).
\end{align*}
Moreover,
\begin{align*}
\widehat{g}(n_1, \dots, n_d)= \int_{\mathbb T^d} g(\boldsymbol{\xi}) \bar{\xi}_1^{n_1} \cdots \bar{\xi}_d^{n_d} d\Theta(\boldsymbol{\xi}) = 0
\end{align*}
unless $n_j \geq 0$ for all $j=1, \cdots, d$ or  $n_j \leq 0$ for all $j=1, \cdots, d$. We call the density $g$ as {\em the phase function} associated to $\varphi$.

Further, the Taylor series expansions of $\log \varphi$ and $H(\cdot, \boldsymbol{\xi})$ around the origin and the following relation
\begin{align*}
\log \varphi (\boldsymbol{z}) = \iota \frac{\pi}{2} + \pi \Im \psi(0) -\iota \pi \int_{\mathbb{T}^d} H(\boldsymbol{z}, \boldsymbol{\xi}) g(\boldsymbol{\xi})d\Theta(\boldsymbol{\xi})
\end{align*}
show that, for any multi-index $\boldsymbol{\alpha} \in \mathbb{N}_0^d$, the Fourier coefficient $\widehat{g}(\boldsymbol{\alpha})$ depends on the Taylor coefficients of $\varphi$ via a universal system of equations
\begin{align} \label{Eq:univ-poly}
\widehat{g}(\boldsymbol{\alpha}) = L_{\boldsymbol{\alpha}} ( c_{\boldsymbol{m}}(\varphi))
\end{align}
where the multi-index $\boldsymbol{m}$ varies over $\Gamma_{\boldsymbol{\alpha}}$. To be precise, $L_{\boldsymbol{\alpha}}$ is a rational function in $|\Gamma_{\boldsymbol{\alpha}}|$ many variables and in the above expression $L_{\boldsymbol{\alpha}} $ is being evaluated at the sequence of Taylor coefficients $(c_{\boldsymbol{m}}(\varphi))_{\boldsymbol{m} \in \Gamma_{\boldsymbol{\alpha}}}$ of $\varphi$ arranged in lexicographical order. We utilize this information in the next section. The exponential representation of a Herglotz-Nevanlinna function was exploited in \cite{Budisic-Putinar} on various domains of $\C^d$, in connection with entropy optimization techniques.

\section{Superresolution phenomenon} \label{Main results}

\subsection{Superresolution in one dimension via Schur's algorithm}
{We start the main part of the present note by an analysis of the single variable setting. Although we shall rely on more sophisticated tools to investigate the superresolution phenomenon in the polydisk, a thorough examination of Schur's algorithm produces a simple version of the superresolution behavior of holomorphic functions depending on one complex variable. To our knowledge the resulting Theorem \ref{Schur} below is new. We start by an auxiliary lemma.}

\begin{lemma}\label{Lemma:superresolution}
Let $f$ be a Blaschke product of degree $n$ with Taylor series around the origin
$$
f(z)=\sum_{j=0}^n a_j z^j + \lito(z^{n+1}).
$$
Then for each $0\leq k \leq n$, the Schur parameter $\gamma_k$ is a smooth function of $2k$ real variables in a neighbourhood of $(\Re a_0, \Im a_0, \dots, \Re a_k, \Im a_k)$. In fact, there is an $\varepsilon >0$ (depending on $a_0, \dots, a_n$) such that in the closed $\varepsilon$-ball, $B_\varepsilon ((\Re a_0, \Im a_0, \dots, \Re a_k, \Im a_k))$ in $\R^{2k}$ with the Euclidean norm, $\gamma_k$ is smooth and its partial derivatives are uniformly bounded.
\end{lemma}
\begin{proof}
Suppose $f$ is a Blaschke product of degree $n$ having Taylor series $$f(z)=\sum_{j=0}^n c_j z^j + \lito(z^{n+1}). $$
 Then the Schur parameters $\{\gamma_k(f) : k\in \mathbb{N}_0 \}$ of $f$ satisfy:
$$
|\gamma_k(f)|<1  \, \text{ for } 0\leq k < n \text{ and } |\gamma_k(f)|=1 \, \text{ for } k\geq n.
$$ 
Note that, the $k^{th}$ Schur parameter only depends on the Taylor coefficients of $f$ up to order $k$, in a smooth manner. Equation \eqref{Eq:Schur parameter} yields
\begin{align} \label{Eq:gamma-k+1}
\gamma_{k+1}(f)= \frac{f_k'(0)}{1-|\gamma_k(f)|^2} \, \text{ for } 0 \leq k <n.
\end{align}
Suppose that the Taylor series around $0$ of $f_k$ is 
\begin{align*}
f_k(z) = \sum_{j=0}^n c^{(k)}_j z^j + \lito(z^{n+1}).
\end{align*}
So, $c_j = c^{(0)}_j$ for $0 \leq j \leq n$. Rewriting \eqref{Eq:Schur parameter} as 
\begin{align*}
zf_{k+1}(1- \overline{\gamma_k(f)} f_k) = f_k - \gamma_k
\end{align*}
and comparing the coefficients of $z^j$ on both sides we find that for $ 1\leq j \leq n-k$,
\begin{align} \label{Eq:c-k-j}
c^{(k)}_j =\left( c^{(k-1)}_{j+1} + \overline{\gamma_{k-1}(f)} \sum_{i=1}^j c^{(k)}_{j-1} c^{(k-1)}_i \right) (1-|\gamma_{k-1}(f)|^2)^{-1} .
\end{align}

Therefore, \eqref{Eq:gamma-k+1} yields 
\begin{align}\label{Eq:gamma-k+2}
\gamma_{k+2}(f)= \frac{c^{(k+1)}_1}{1-|\gamma_{k+1}(f)|^2} = \frac{c^{(k)}_2 + \overline{\gamma_{k}(f)} \gamma_{k+1}(f) c^{(k)}_1 }{(1-|\gamma_{k}(f)|^2)(1-|\gamma_{k+1}(f)|^2)}
\end{align}
for $0\leq k \leq n-2$. 

Let $c_k= x_k + \iota y_k$ for $0\leq k \leq n$ where $x_k, y_k \in \R$. Clearly, $\gamma_0(f)=c_0$, $\gamma_1(f)= c_1/(1-|c_0|^2)$ and hence both are rational functions in $x_0, y_0, x_1, y_1$ as variables. Therefore, from \eqref{Eq:gamma-k+2} we observe that $\gamma_2(f)$ is a rational function in $x_0, y_0, x_1, y_1, x_2, y_2$ as variables. Summing up, we proved inductively that $\gamma_{k+2}(f)$ is a rational function in $x_0, y_0, \dots, x_{k+2}, y_{k+2}$ as variables, where $0\leq k \leq n-2$. 

Let us denote these rational functions by $\psi_0, \dots, \psi_n$ where
$$
\psi_k(x_0, y_0, \dots, x_{k}, y_{k})= \gamma_k(f)
$$
for $0\leq k \leq n$. 

Since $|\gamma_j(f)|<1$ for $0\leq j <n$, the functions $\psi_k$ do not have singularity at $(\Re a_0, \Im a_0, \dots, \Re a_k, \Im a_k)$ for each $k$. Therefore, we can choose an $\varepsilon >0$ such that for each $0 \leq k \leq n$, $\psi_k$ is smooth on the $\varepsilon$-ball,
$B_\varepsilon((\Re a_0, \Im a_0, \dots, \Re a_k, \Im a_k))$ in $\R^{2k}$. Consequently, their partial derivatives are uniformly bounded on $B_\varepsilon((\Re a_0, \Im a_0, \dots, \Re a_k, \Im a_k))$.

\end{proof}

The following theorem describes superresolution behaviour of Schur functions near a finite Blaschke product.
\begin{theorem}\label{Schur}
Let $f$ be a Blaschke product of degree $n$ with Taylor series around the origin
\begin{align}\label{Eq:Taylor-series f}
f(z)=\sum_{j=0}^n c_j z^j + \lito(z^{n+1}).
\end{align}
 There exist three positive real numbers $\varepsilon$, $L$ and $M$ (depending on the Taylor coefficients of $f$ up to order $n$) such that for every $g \in \mathcal{S}(\D)$ with Taylor series around the origin
\begin{align} \label{Eq:Taylor-series g}
g(z)=\sum_{j=0}^n \tilde{c}_j z^j + \lito(z^{n+1})
\end{align}
 and $\sum_{j=0}^n |c_j - \tilde{c}_j|^2  < \varepsilon^2$, we have  
\begin{align*}
|f(z)-g(z)| \leq \frac{4 M \varepsilon}{(L-M\varepsilon)(1-|z|)} \, \, \text{ for every } z\in \D.
\end{align*}
\end{theorem}

\begin{proof}
For the Blaschke product $f$ as above, applying Lemma~ \ref{Lemma:superresolution} we obtain an $\varepsilon >0$ (depending on the Taylor coefficients of $f$ up to order $n$) such that in the $\varepsilon$-ball, $B_\varepsilon((\Re c_0, \Im c_0, \dots, \Re c_k, \Im c_k))$ the Schur parameter $\gamma_k$ is smooth with uniformly bounded partial derivatives. 

Suppose $g \in \mathcal{S}(\D)$ with Taylor series as in \eqref{Eq:Taylor-series g} and $\sum_{j=0}^n |c_j - \tilde{c}_j|^2  < \varepsilon^2$. We denote the Wall polynomials corresponding to $f$ and $g$ by $A_k, B_k$ and $\tilde{A}_k, \tilde{B}_k$, respectively. By Proposition \ref{Prop:Wall poly}, there exists a $h \in \mathcal{S}(\D)$ such that 
\begin{align} \label{Eq:Wall-convergent}
g = \frac{\tilde{A}_n + z \tilde{B}_n^* h}{\tilde{B}_n + z\tilde{A}_n^* h}. 
\end{align}
Also, $f$ being a Blaschke product, $f= \frac{A_n}{B_n}$. Recall that, $A_k$ and $B_k$ are polynomials of degree $k$ with coefficients depending on the Schur parameters $\gamma_0(f), \dots, \gamma_k(f)$. Further,
$\gamma_k(f) = \psi_k(\Re c_0, \Im c_0, \dots, \Re c_k, \Im c_k)$ for  $0\leq k \leq n$.  Thus the coefficients of the Wall polynomials of degree $k$ corresponding to a Schur function having $n^{th}$ degree Taylor polynomial $\sum_{j=0}^n a_j z^j$ are smooth functions of $\Re a_0, \Im a_0, \dots, \Re a_k, \Im a_k$ in $B_\varepsilon((\Re c_0, \Im c_0, \dots, \Re c_k, \Im c_k))$ in $\R^{2k}$. Since the partial derivatives of the functions $\psi_k$ are bounded on $B_\varepsilon((\Re c_0, \Im c_0, \dots, \Re c_k, \Im c_k))$, the coefficients of the Wall polynomials are Lipschitz continuous.

Therefore, there exists $M>0$ (chosen to be maximum of all such Lipschitz constants) depending on $c_0, \dots, c_n$ such that 
\begin{align}\label{Eq:Wall-estimate}
 \|A_n - \tilde{A}_n \|_\infty \leq M \varepsilon \,\, \text{ and } \|B_n - \tilde{B}_n\|_\infty \leq M \varepsilon
\end{align}
 as $\sum_{j=0}^n |c_j - \tilde{c}_j|^2  < \varepsilon^2$.
 
 Again, the modulus of a polynomial $p$ and its reflection $p^*$ are the same on $\mathbb{T}$,  hence $\|p\|_\infty = \|p^*\|_\infty.$ Thus, \eqref{Eq:Wall-estimate} implies
 \begin{align} \label{Eq:Wall-estimate-refl}
 \|A_n^* - \tilde{A}_n^* \|_\infty \leq M \varepsilon \,\, \text{ and } \|B_n^* - \tilde{B}_n^*\|_\infty \leq M \varepsilon.
\end{align}
Now, for $z\in \D$,
\begin{align}\label{Eq:Main-estimate}
&|f(z)-g(z)| \notag \\ 
= &\left| \frac{A_n}{B_n} - \frac{\tilde{A}_n + z \tilde{B}_n^* h}{\tilde{B}_n + z\tilde{A}_n^* h} \right| \,\, \text{ (by \eqref{Eq:Wall-convergent})}  \notag \\
= &\left|   \frac{A_n}{B_n} \frac{\tilde{B}_n - B_n}{\tilde{B}_n + z\tilde{A}_n^* h}  +  \frac{ A_n - \tilde{A}_n }{\tilde{B}_n + z\tilde{A}_n^* h} + zh  \frac{A_n}{B_n}  \frac{\tilde{A}_n^* - A_n^*}{\tilde{B}_n + z\tilde{A}_n^* h} \right. \notag\\
 &\left. + zh \frac{A_n A_n^* - B_n B_n^*}{B_n(\tilde{B}_n + z\tilde{A}_n^* h)}  + zh \frac{B_n^* - \tilde{B}_n^*}{\tilde{B}_n + z\tilde{A}_n^* h} \right| 
\end{align}
 
 Since $B_n$ is non-vanishing in $\overline{\D}$, $L= \inf \{|B_n(z)|: z\in \overline{\D} \} >0$. Now using \eqref{Eq:Wall-estimate}, we have
 \begin{align} \label{Eq:lower-estimate-1}
 |\tilde{B}_n(z)|> |B_n(z)|- M\varepsilon \geq L - M\varepsilon.
 \end{align}
 Again, \eqref{Eq:lower-estimate-1} and part ($3$) of Proposition \ref{Prop:Wall poly} yield
 \begin{align}\label{Eq:lower-estimate-2}
 |\tilde{B}_n + z\tilde{A}_n^* h| &  \geq |\tilde{B}_n| \left( 1- \left| zh \frac{\tilde{A}_n^*}{\tilde{B}_n}\right| \right)  \notag \\
                                  &   \geq (L - M\varepsilon) (1-|z|).
 \end{align}

Further, part ($1$) of Proposition \ref{Prop:Wall poly} implies that $B_n^* B_n - A_n^* A_n=0$ as $|\gamma_n(f)|=1$. This fact along with the estimates in 
\eqref{Eq:Main-estimate}, \eqref{Eq:Wall-estimate} and \eqref{Eq:Wall-estimate-refl} finally prove 
 \begin{align*}
 |f(z)-g(z)| \leq \frac{4M\varepsilon}{(L - M\varepsilon) (1-|z|)}
 \end{align*}
 as required. 
\end{proof}

\subsection{Superresolution in the polydisk}
We turn now to the main topic of this article: the superresolution phenomenon for Herglotz-Nevanlinna functions in the polydisc $\D^d$. Two main ingredients enter into the analysis below:
 the exponential representation of a Herglotz-Nevanlinna function and volume estimates of principal semi-algebraic sets. Before stating the main result we fix some notations. 

For a holomorphic function $f$ on $\D^d$, the $\boldsymbol{\alpha}$-th Taylor coefficient (at the origin) is denoted by $c_{\boldsymbol{\alpha}}(f)$ and its Taylor polynomial of multi-degree at most $\boldsymbol{n}$ is denoted by $T_{\boldsymbol{n}}(f)$ i.e., 
$$
T_{\boldsymbol{n}}(f) = \sum_{\boldsymbol{\alpha} \in \Gamma_{\boldsymbol{n}} } c_{\boldsymbol{\alpha}} (f) \boldsymbol{z}^{\boldsymbol{\alpha}}.
$$
Identifying the Taylor polynomial, $T_{\boldsymbol{n}}(f)$ with the vector $(c_{\boldsymbol{\alpha}}(f))_{\boldsymbol{\alpha} \in \Gamma_{\boldsymbol{n}}}$ in $\mathbb{C}^{|\Gamma_{\boldsymbol{n}}|}$, we define the norm
$$
\|T_{\boldsymbol{n}}(f)\|_2 = \sum_{\boldsymbol{\alpha} \in \Gamma_{\boldsymbol{n}} } |c_{\boldsymbol{\alpha}} (f)|^2.
$$
Now we are ready to state our main result.
\begin{theorem}\label{Thm:supres-polydisk}
Let $R$ be a Cayley rational inner function of multi-degree $\boldsymbol{n}= (n_1, \dots, n_d)$ on $\D^d$ and let $f$ be in $\mathcal{H}(\D^d)$. Then for each compact subset $K$ of $\D^d$, there exist constants $C(K, R)>0$ (depending on $K$ and $R$), $\rho >0$ (depending on $R$) and a positive exponent $\kappa$ (depending on $R$), such that 

\begin{align*}
\|R - f\|_{\infty, K}^\kappa \ \leq \ C(K, R) \|T_{\boldsymbol{n}}(R) - T_{\boldsymbol{n}}(f) \|_2 
\end{align*}
whenever $\|T_{\boldsymbol{n}}(R) - T_{\boldsymbol{n}}(f) \|_2 < \rho $.

\end{theorem}
The proof is split into several lemmas of potential independent interest.

\begin{lemma}\label{Lemma:1}
Let $R$ be as in Theorem \ref{Thm:supres-polydisk} and let $g \in L^1(\mathbb{T}^d), \ 0\leq g \leq 1$ be the phase function of $1/2 + (\iota/\pi) \log R $. Then there exists a trigonometric polynomial $P$ on $\mathbb{T}^d$ of multi-degree at most $\boldsymbol{n}$ such that
\begin{align*}
g = \chi_{S} \ \text{ with } S = \{\boldsymbol{\xi} \in \mathbb{T}^d: P(\boldsymbol{\xi}, \overline{\boldsymbol{\xi}}) >0\}, 
\end{align*}
\end{lemma}

Above $\chi_S$ denotes the characteristic function of $S$.

\begin{proof}
We already noted in Subsection \ref{subsection:Finite-determinateness} that any Cayley rational inner function of multi-degree $\boldsymbol{n}$ and hence $R$ is completely determined by its Taylor polynomial $T_{\boldsymbol{n}} (R)$ of multi-degree $\boldsymbol{n}$. Now, as in Subsection \ref{subsection:Exp-transform} we look at $\psi = \frac{1}{2} + \frac{\iota}{\pi} \log R $ and the phase function $g$ in $L^1(\mathbb{T}^d)$ with $0\leq g \leq 1$ such that
\begin{align} \label{Eq:exp-transform}
R(\boldsymbol{z})= \iota \exp ( \pi \Im \psi(0))  \exp \left(-\iota \pi \int_{\mathbb{T}^d} H(\boldsymbol{z}, \boldsymbol{\xi}) g(\boldsymbol{\xi})d\Theta(\boldsymbol{\xi}) \right)
\end{align}
 and  
\begin{align*}
\Re \psi(\boldsymbol{z}) = \int_{\mathbb{T}^d} \mathcal{P}_{\D^d}(\boldsymbol{z}, \boldsymbol{\xi}) g(\boldsymbol{\xi})d\Theta(\boldsymbol{\xi}).
\end{align*}
Let $$\widehat{g}(\boldsymbol{\alpha}) = \int_{\mathbb{T}^d} g(\boldsymbol{\xi}) \bar{\boldsymbol{\xi}}^{\boldsymbol{\alpha}} d\Theta(\boldsymbol{\xi})$$ for $\boldsymbol{\alpha}$ in $\mathbb{Z}^d$. But we have from \eqref{Eq:Pluri-measure}, $\widehat{g}(\boldsymbol{\alpha}) = 0 $ if $\boldsymbol{\alpha} \in \mathbb{Z}^d \setminus (\mathbb{N}_0^d \cup - \mathbb{N}_0^d )$.

Thus for $\boldsymbol{z} \in \mathbb{D}^d$ we have
\begin{align*}
\Re \psi(\boldsymbol{z}) & =  \sum_{\boldsymbol{\alpha} \in \mathbb{N}_0^d} \widehat{g}(\boldsymbol{\alpha})  \boldsymbol{z}^{\boldsymbol{\alpha}} + \overline{\widehat{g}(\boldsymbol{\alpha})} \ \overline{\boldsymbol{z}}^{\boldsymbol{\alpha}}.
\end{align*} 

Again, \eqref{Eq:exp-transform} shows that the Fourier coefficients $\{\widehat{g}(\boldsymbol{\alpha}):  \boldsymbol{\alpha} \in \Gamma_{\boldsymbol{n}} \}$ determine $T_{\boldsymbol{n}} (R)$ and hence $R$ completely. Thus  
$$
\{\widehat{g}(\boldsymbol{\alpha}):  \boldsymbol{\alpha} \in \Gamma_{\boldsymbol{n}} \}
$$
entirely determines $\psi$ and hence $g$.
Therefore, the density $g$ is finitely determined by the real moments 
\begin{align*}
\Re \widehat{g}(\boldsymbol{\alpha}) =  \int_{\mathbb{T}^d} g(\boldsymbol{\xi}) \frac{(\overline{\boldsymbol{\xi}}^{\boldsymbol{\alpha}}+ \boldsymbol{\xi}^{\boldsymbol{\alpha}})}{2} d\Theta(\boldsymbol{\xi}),
\end{align*}
and 
\begin{align*}
\Im \widehat{g}(\boldsymbol{\alpha}) =  \int_{\mathbb{T}^d} g(\boldsymbol{\xi})  \frac{ (\overline{\boldsymbol{\xi}}^{\boldsymbol{\alpha}}- \boldsymbol{\xi}^{\boldsymbol{\alpha}})}{2\iota} d\Theta(\boldsymbol{\xi}),
\end{align*}
for $ \boldsymbol{\alpha} \in \Gamma_{\boldsymbol{n}}$.

For any $\eta \in L^1(\mathbb{T}^d, d\Theta)$ define, 
$$
\boldsymbol{a}(\eta) = (\Re a_{\boldsymbol{\alpha}}, \Im a_{\boldsymbol{\alpha}})_{\boldsymbol{\alpha} \in \Gamma_{\boldsymbol{n}}} \in \mathbb{R}^{2|\Gamma_{\boldsymbol{n}}|}, \ \text{where } a_{\boldsymbol{\alpha}} = \widehat{\eta}(\boldsymbol{\alpha}).
$$

We consider the closed, convex and bounded subset of $\mathbb{R}^{2|\Gamma_{\boldsymbol{n}}|}$
\begin{align*}
\Sigma_{\boldsymbol{n}}= \left \lbrace \boldsymbol{a}(\eta) : \eta \in L^1(\mathbb{T}^d, d\Theta), \ 0 \leq \eta \leq 1 \right \rbrace.
\end{align*}
We view every linear functional $\mathcal{F}$ on $\mathbb{R}^{2|\Gamma_{\boldsymbol{n}}|}$ as a trigonometric polynomial 
$$
J(\boldsymbol{\xi}, \overline{\boldsymbol{\xi}}) = \sum_{\boldsymbol{\alpha} \in \Gamma_{\boldsymbol{n}}} b_{\boldsymbol{\alpha}} \frac{(\overline{\boldsymbol{\xi}}^{\boldsymbol{\alpha}}+ \boldsymbol{\xi}^{\boldsymbol{\alpha}})}{2} +  
\tilde{b}_{\boldsymbol{\alpha}}  \frac{ (\overline{\boldsymbol{\xi}}^{\boldsymbol{\alpha}}- \boldsymbol{\xi}^{\boldsymbol{\alpha}})}{2\iota}
$$
on $\mathbb{T}^d$ such that 
$$
\mathcal{F}(\boldsymbol{a}(\eta)) = \int_{\mathbb{T}^d} J \eta d\Theta = \sum_{\boldsymbol{\alpha} \in \Gamma_{\boldsymbol{n}}} b_{\boldsymbol{\alpha}} \Re a_{\boldsymbol{\alpha}} + \tilde{b}_{\boldsymbol{\alpha}} \Im a_{\boldsymbol{\alpha}}.
$$

 Returning to the phase function $g$, we observe that the sequence of moments $\boldsymbol{a}(g)$ in $\Sigma_{\boldsymbol{n}}$ is a boundary point of $\Sigma_{\boldsymbol{n}}$. Indeed, if $\boldsymbol{a}(g)$ is an interior point of $\Sigma_{\boldsymbol{n}}$ then there is a linear functional $\mathcal{F}$ on $\mathbb{R}^{2|\Gamma_{\boldsymbol{n}}|}$ as above such that
 $$
 \mathcal{F}(\boldsymbol{a}(g)) =  \int_{\mathbb{T}^d} J g d\Theta <  \| J\|_1.
 $$
 This strict inequality contradicts the fact that $g$ is determined by the moments $\boldsymbol{a}(g)$. Thus
there exists a unique trigonometric polynomial
\begin{align} \label{Eq:key-pol}
P(\boldsymbol{\xi}, \overline{\boldsymbol{\xi}}) = \sum_{\boldsymbol{\alpha} \in \Lambda_{\boldsymbol{n}}} r_{\boldsymbol{\alpha}} \frac{(\overline{\boldsymbol{\xi}}^{\boldsymbol{\alpha}}+ \boldsymbol{\xi}^{\boldsymbol{\alpha}})}{2} +  
\tilde{r}_{\boldsymbol{\alpha}}  \frac{ (\overline{\boldsymbol{\xi}}^{\boldsymbol{\alpha}}- \boldsymbol{\xi}^{\boldsymbol{\alpha}})}{2\iota}
\end{align} 
with $r_{\boldsymbol{\alpha}}, \tilde{r}_{\boldsymbol{\alpha}}  \in \mathbb{R}$ such that 
$$
\int_{\mathbb{T}^d} P g d\Theta =  \| P\|_1.
$$
An analysis of the equality sign in the abstract $L^1-L^\infty$ duality implies that
\begin{align*}
g = \chi_{\{P >0\} } \ \ \Theta- \text{a.e. on } \mathbb{T}^d.
\end{align*}
Here $\chi_{ \{P >0\} }$ denotes the characteristic function of the set $$ S= \{\boldsymbol{\xi} \in \mathbb{T}^d: P(\boldsymbol{\xi}, \overline{\boldsymbol{\xi}}) >0\}.$$
For a discussion of the rigidity of "black and white" shades in the $L$-problem of moments see \cite[Chapter VII]{Krein-Nudelman}. The specific instance of semi-algebraic sets in Euclidean space is treated in \cite{Putinar}. \

This completes the proof of Lemma \ref{Lemma:1}.
\end{proof}

\subsection*{Decomposing $\mathbb T^d$}
We parametrize the unit circle $\mathbb{T}$ by the following two coordinate charts:
\begin{align*}
\Psi_1 : [-1,1] \rightarrow \mathbb{T}_1 \bydef \{e^{\iota \theta}: \theta \in [-\pi/2, \pi/2] \} \subseteq \mathbb{T}, \ \text{ where }
\end{align*}
$$ \Psi_1(t) = \frac{1-t^2}{1+t^2} + \iota \frac{2t}{1+t^2}, $$
and 
\begin{align*}
\Psi_2 : [-1,1] \rightarrow \mathbb{T}_2 \bydef \{e^{\iota \theta}: \theta \in [\pi/2, 3\pi/2] \} \subseteq \mathbb{T}, \ \text{ where }
\end{align*}
\begin{align*}
\Psi_2(t) \bydef - \frac{1-t^2}{1+t^2} - \iota \frac{2t}{1+t^2} \text{ on } [-1,1].
\end{align*}

Note that, $ \Psi_j ([-1, 1]) = \mathbb{T}_j .$ Also, $\mathbb{T} = \mathbb{T}_1 \cup \mathbb{T}_2$. Therefore 
\begin{align}\label{Eq:torus-decom}
\mathbb{T}^d =  \cup_{X_1, \dots, X_d \in \{\mathbb{T}_k: k=1,2\}} \ (X_1 \times \dots \times X_d ).
\end{align}
For simplicity we rewrite the above union as  
$$
\mathbb{T}^d =  \cup_{j=1}^{2^d} S_j
$$
where the sets, $S_j$ are the pieces of $\mathbb{T}^d$ appearing in \eqref{Eq:torus-decom}. 
Consequently each $S_j$ has a parametrization from $[-1, 1]^d$ through the maps $\Psi_j$ where each component of $S_j$ is parametrized by either $\Psi_1$ or $\Psi_2$ depending on the component. It is easy to see that the $\Theta$ measure of $S_j \cap S_k \subseteq \mathbb{T}^d $ is zero as $\mathbb{T}_1 \cap \mathbb{T}_2 = \{\iota, -\iota\}$. Also, $ S =  \cup_{j=1}^{2^d} (S \cap S_j)$. 

In view of the decomposition above, the phase function $g$ and the set $S$ appearing in Lemma \ref{Lemma:1} can be expressed as:
\begin{align*}
g= \chi_S = \sum_{j=1}^{2^d} \chi_{S \cap S_j} \ \ \Theta-\text{a.e.}.
\end{align*}
Let $\lambda$ denote the Lebesgue measure supported on $\Delta \bydef [-1, 1]^d$. Now, our next lemma is the following.

\begin{lemma}\label{Lemma:2}
Fix $1\leq j \leq 2^d$. Suppose $\Phi = (\Psi_{j_1}, \dots, \Psi_{j_d})$ be the parametrization of $S_j$ from $\Delta$ where each $j_k \in \{1,2\}$. Let $g_j \bydef  \chi_{S \cap S_j}$ and $\tilde{g} \bydef g_j \circ \Phi $. Then, for every $u \in L^1(\mathbb{T}^d \cap S_j)$,
\begin{align*}
\| g_j - u \|_{L^1(\mathbb{T}^d \cap S_j)} \ \leq \ 2^d \| \tilde{g}_j - u \circ \Phi \|_{L^1(\lambda)}.
\end{align*}
\end{lemma}

\begin{proof}
For $\boldsymbol{t} \in \Delta$, $\tilde{g}(\boldsymbol{t})  = \chi_{S \cap S_j } \circ \Phi(\boldsymbol{t})$ and so,
                    \[ 
                          \tilde{g}(\boldsymbol{t}) = 
                  \begin{cases}
                        1, \ \text{ if } P(\Phi(\boldsymbol{t}))>0 \\
                         0, \ \text{ otherwise.} \\ 
                  \end{cases}
\]

From \eqref{Eq:key-pol}, we have 
\begin{align*}
 P(\Phi(\boldsymbol{t})) =  \sum_{\boldsymbol{\alpha} \in \Lambda_{\boldsymbol{n}}} r_{\boldsymbol{\alpha}} \frac{(\overline{\Phi(\boldsymbol{t})}^{\boldsymbol{\alpha}}+ \Phi(\boldsymbol{t})^{\boldsymbol{\alpha}})}{2} +  
\tilde{r}_{\boldsymbol{\alpha}}  \frac{ (\overline{\Phi(\boldsymbol{t})}^{\boldsymbol{\alpha}}- \Phi(\boldsymbol{t})^{\boldsymbol{\alpha}})}{2\iota}
\end{align*}
which can be rewritten as 
\begin{align}\label{Eq:P-expression}
P(\Phi(\boldsymbol{t})) = \frac{Q(\boldsymbol{t})}{\prod_{j=1}^d (1+t_j^2)^{|\boldsymbol{n}|}} 
\end{align}
where $Q$ is a real polynomial in $\mathbb{R}[\boldsymbol{t}]$ of total degree atmost $2d |\boldsymbol{n}| $. This implies that $P(\Phi(\boldsymbol{t})) >0 $ if and only if $Q(\boldsymbol{t})>0$ for $\boldsymbol{t} \in \Delta $ and hence 
\begin{align*}
\tilde{g}_j= \chi_{\{ Q>0\}} \ \  \mu-\text{a.e. on }  \Delta.
\end{align*}
Here $\mu$ is the finite measure supported on $[-1, 1]^d$ given by $d\mu(\boldsymbol{t}) = 2^d \frac{d\lambda(\boldsymbol{t})}{(1+t_1^2) \dots (1+t_d^2)}$. Note that, $\mu$ is the pull-back of the measure $\Theta$ on $\mathbb{T}^d$ under the parametrization $\Phi$. Also, the measures $\lambda$ and $\mu$ are absolutely continuous with respect each other on $\Delta$. Thus 
\begin{align*}
\tilde{g}_j= \chi_{\{ Q>0\}} \ \  \lambda-\text{a.e. on }  \Delta
\end{align*}
as well. Moreover,
\begin{align}
\| g_j - u \|_{L^1(\mathbb{T}^d \cap S_j)} & = \int_{\mathbb{T}^d \cap S_j} |g_j \circ \Phi(\boldsymbol{t}) - u \circ \Phi(\boldsymbol{t})| d\mu(\boldsymbol{t}) \notag \\
& \leq 2^d \int_{\Delta} |\tilde{g}_j(\boldsymbol{t}) - u \circ \Phi(\boldsymbol{t})| d\lambda(\boldsymbol{t}) \notag \\
& = 2^d \| \tilde{g}_j - u \circ \Phi \|_{L^1(\lambda)}
\end{align}
for every $u \in L^1(\mathbb{T}^d \cap S_j)$.
\end{proof}

For the fixed $j$, $\tilde{g}_j$ and the real polynomial $Q$ as in Lemma \ref{Lemma:2}, we denote the coefficient of $\boldsymbol{t}^{\boldsymbol{\beta}}$ in $Q$ by $q_{\boldsymbol{\beta}}$ for every multi-index $\boldsymbol{\beta}$. Let $\boldsymbol{m}=(m_1, \dots, m_d)$ be an {\em admissible index} for the polynomial $Q$ in the sense that $q_{\boldsymbol{m}}$ is non-zero and there exists a permutation $(\sigma(1), \dots, \sigma(d))$ of $(1, \dots, d)$ such that for every multi-index $\boldsymbol{\beta}=(\beta_1, \dots, \beta_d)$ with $q_{\boldsymbol{\beta}} \neq 0$, either $m_{\sigma(1)} > \beta_{\sigma(1)}$, or there exists an index $l$, $l\geq 2$, satisfying $m_{\sigma(j)} > \beta_{\sigma(j)}$ and $m_{\sigma(k)} = \beta_{\sigma(k)}$ for $1\leq k \leq l-1$. See \cite[Section 3]{Putinar} for more details. Then we have the following estimation.

\begin{lemma}[Volume estimate] \label{Lemma:3}
For $j$, $\tilde{g}_j$, $Q$ and $\boldsymbol{m}$ as above, there exist a constant $C_j$ (depending only on $d$) 
such that  
\begin{align}\label{Eq:Lewis-estimate}
\int_{\Delta} (\tilde{g}_j - v)Q d\lambda \geq \frac{\|\tilde{g}_j - v \|_{L^1(\lambda)}^{|\boldsymbol{m}|+1}}{C_j^{|\boldsymbol{m}|}(1+ |\boldsymbol{m}|)},
\end{align}
whenever $v \in L^1(\lambda)$ with $0 \leq v \leq 1$ and
$
\|\tilde{g}_j - v \|_{L^1(\lambda)} \leq C_j \frac{|q_{\boldsymbol{m}}|^{1/|\boldsymbol{m}|}}{(8 |\boldsymbol{n}|d)}.
$
\end{lemma}

\begin{proof}
 Following \cite{Lewis, Putinar}, define
$$
\Lambda_{|Q|}(\varepsilon) = \inf \left \lbrace  \int_{\Delta} |Q| v d\lambda : v\in L^1(\lambda), 0 \leq v \leq 1, \int_{\Delta} v d\lambda \geq \varepsilon   \right \rbrace.
$$
For every $v \in L^1(\lambda)$ with $0 \leq v \leq 1$,  
\begin{align} \label{Eq:estimate-convergence}
\int_{\Delta}  \prod_{j=1}^d (1+t_j^2)^{|\boldsymbol{n}|} (\tilde{g}_j - v) (P \circ \Phi) d\lambda = \int_{\Delta} (\tilde{g}_j - v)Q d\lambda \geq \Lambda_{|Q|}(\|\tilde{g}_j - v \|_{L^1(\lambda)}). 
\end{align}
The first equality is obtained using \eqref{Eq:P-expression} and the for the elementary proof of the second inequality we refer to \cite[Lemma~2.10]{Lewis}. 

Now, proceeding with the argument as in Section~4 of \cite{Putinar}, we obtain a constant $C_j$ (depending only on $d$) such that for every $v \in L^1(\lambda)$ with $0 \leq v \leq 1$,  
\begin{align*}
\int_{\Delta} (\tilde{g}_j - v)Q d\lambda \geq \frac{\|\tilde{g}_j - v \|_{L^1(\lambda)}^{|\boldsymbol{m}|+1}}{C_j^{|\boldsymbol{m}|}(1+ |\boldsymbol{m}|)},
\end{align*}
whenever 
$$
\|\tilde{g}_j - v \|_{L^1(\lambda)} \leq C_j \frac{|q_{\boldsymbol{m}}|^{1/|\boldsymbol{m}|}}{(8 |\boldsymbol{n}|d)},
$$
completing a sketch of the proof of Lemma \ref{Lemma:3}. We refer to \cite{Putinar} for full details.
\end{proof}

\noindent \textbf{Note:} Recall that, so far we have worked with $1\leq j \leq 2^d$ which has been fixed in Lemma \ref{Lemma:2} and obtained the polynomial $Q$ as well as the corresponding admissible index $\boldsymbol{m}$ by doing the analysis on the piece $S_j \subseteq \mathbb{T}^d$. In order to keep track of this dependence on $j$ we rename the polynomial $Q$ and the index $\boldsymbol{m}$ as $Q_j$ and $\boldsymbol{m}_j$, respectively. Also, the coefficient of $\boldsymbol{t}^{\boldsymbol{\beta}}$ of $Q_j$ is relabelled as $q_{\boldsymbol{\beta}}^{(j)}$. Further, we rename the parametrization map, $\Phi$ on $S_j$ as $\Phi_j$.  

Thus the estimate in \eqref{Eq:Lewis-estimate} of Lemma \ref{Lemma:3} with the help of \eqref{Eq:P-expression} can be rewritten as
\begin{align}\label{Eq:vol-estimate}
&\int_{\Delta} (\tilde{g}_j - v)Q_j d\lambda \ \geq \ \frac{\|\tilde{g}_j - v \|_{L^1(\lambda)}^{|\boldsymbol{m}_j|+1}}{C_j^{|\boldsymbol{m}_j|}(1+ |\boldsymbol{m}_j|)} \notag \\
\text{i.e.,} & \int_{\Delta}  \prod_{j=1}^d (1+t_j^2)^{|\boldsymbol{n}|} (\tilde{g}_j - v) (P_j \circ \Phi_j) d\lambda \ \geq \ \frac{\|\tilde{g}_j - v \|_{L^1(\lambda)}^{|\boldsymbol{m}_j|+1}}{C_j^{|\boldsymbol{m}_j|}(1+ |\boldsymbol{m}_j|)},
\end{align} 
for every $v \in L^1(\lambda)$ with $0 \leq v \leq 1$ such that 
$$
\|\tilde{g}_j - v \|_{L^1(\lambda)} \leq C_j \frac{|q_{\boldsymbol{m}}^{(j)}|^{1/|\boldsymbol{m}_j|}}{8d |\boldsymbol{n}|}.
$$

Now we are in a position to complete the proof of the superresolution as stated in Theorem \ref{Thm:supres-polydisk}.
\begin{proof}[\textbf{Proof of Theorem \ref{Thm:supres-polydisk}}]
Consider a compact subset $K$ of $\D^d$ and a Herglotz-Nevanlinna function $f$ on $\D^d$. Let $h$ be the phase function corresponding to the Herglotz-Nevanlinna function $\varphi = \frac{1}{2} + \frac{\iota}{\pi} \log f$. Using the exponential transform, we have $h\in L^1(\mathbb{T}^d)$ with $0 \leq h \leq 1$ such that
\begin{align} \label{Eq:exp-f}
f(\boldsymbol{z})= \iota \exp ( \pi \Im \varphi(0))  \exp \left(-\iota \pi \int_{\mathbb{T}^d} H(\boldsymbol{z}, \boldsymbol{\xi}) h(\boldsymbol{\xi})d\Theta(\boldsymbol{\xi}) \right),
\end{align}
for $\boldsymbol{z} \in \D^d$. Also, in Lemma \ref{Lemma:1} we have the phase function $g$ of $\psi$ such that \eqref{Eq:exp-transform} holds for the Cayley rational inner function $R$.
For $\boldsymbol{z} \in K$ and $\boldsymbol{\xi} \in \mathbb{T}^d$,
$|H(\boldsymbol{z}, \boldsymbol{\xi} )| \leq \frac{2}{\prod_{l=1}^d (1-|z_l|)} + 1 \leq C'(K)$, for some positive constant $C'(K)$ depending on the compact set $K$.

Therefore, for $\boldsymbol{z} \in K$, \eqref{Eq:exp-transform} and \eqref{Eq:exp-f} imply
\begin{align}\label{Eq:pointwise-estimate}
& |R(\boldsymbol{z}) - f(\boldsymbol{z})| \notag \\
& = \left|  \exp ( \pi \Im \psi(0))  \exp \left(-\iota \pi \int_{\mathbb{T}^d} H(\boldsymbol{z}, \boldsymbol{\xi}) g(\boldsymbol{\xi})d\Theta(\boldsymbol{\xi}) \right) \right. \notag \\
& \left. - \exp ( \pi \Im \varphi(0))  \exp \left(-\iota \pi \int_{\mathbb{T}^d} H(\boldsymbol{z}, \boldsymbol{\xi}) h(\boldsymbol{\xi})d\Theta(\boldsymbol{\xi}) \right)  \right| \notag \\
& \leq C(K)  \left| \pi \Im (\psi(0) - \varphi(0)) + \iota \pi \int_{\mathbb{T}^d} H(\boldsymbol{z}, \boldsymbol{\xi}) (h-g)(\boldsymbol{\xi}) d\Theta(\boldsymbol{\xi})  \right|    \notag \\
&  \leq \pi C(K) |\psi(0) - \varphi(0)| +   \pi C'(K) \| g -h \|_{L^1(\mathbb{T}^d)},
\end{align}
where the positive constants $C(K)$ is obtained by the mean value theorem for the exponential function on the compact set $K$.
Again, the almost disjoint decomposition of $\mathbb{T}^d$ using the sets $\{S_j: 1\leq j \leq 2^d \}$ implies that
\begin{align}\label{Eq:Mainak-2}
\| g -h \|_{L^1(\mathbb{T}^d)}  = \sum_{j=1}^{2^d} \| g_j -h \chi_{S_j} \|_{L^1(\mathbb{T}^d \cap S_j)} \leq 2^d \sum_{j=1}^{2^d}  \| \tilde{g}_j - (h \chi_{S_j}) \circ \Phi_j \|_{L^1(\lambda)}.
\end{align}

Let us denote $h \chi_{S_j} \circ \Phi_j$ by $v_j$. The above equality implies
\begin{align}\label{Eq:Putinar-estimate}
  & \| g -h \|_{L^1(\mathbb{T}^d)}  \notag \\
  & \leq \  2^d \sum_{j=1}^{2^d}  \| \tilde{g}_j - v_j \|_{L^1(\lambda)} ^{|\boldsymbol{m}_j|+1} \notag \\
  & \leq \  2^d   \sum_{j=1}^{2^d} \left[ C_j^{|\boldsymbol{m}_j|}(1+ |\boldsymbol{m}_j|) \int_{\Delta}  \prod_{j=1}^d (1+t_j^2)^{|\boldsymbol{n}|} (\tilde{g}_j - v) (P_j \circ \Phi_j) \ d\lambda \right]^{\frac{1}{|\boldsymbol{m}_j|+1}} \ \ (\text{using } \eqref{Eq:vol-estimate}) \notag \\
  & \leq C \sum_{j=1}^{2^d} \left[ \int_{\Delta} (\tilde{g}_j - v) (P_j \circ \Phi_j)\  d\lambda \right]^{\frac{1}{|\boldsymbol{m}_j|+1}}  \notag \\
  & \leq C' \left[ \sum_{j=1}^{2^d}  \int_{\Delta} (\tilde{g}_j - v) (P_j \circ \Phi_j) \ d\lambda \right]^{\frac{1}{\kappa}}  \notag \\
  & \leq C' \left[  \sum_{j=1}^{2^d}  \int_{\mathbb{T}^d \cap S_j} (\tilde{g}_j \circ \Phi_j^{-1} - v \circ \Phi_j^{-1}) P_j \ d\Theta \right]^{\frac{1}{\kappa}} \ \ \text{(by change of variables)}  \notag \\
  & = C' \left[  \sum_{j=1}^{2^d}  \int_{\mathbb{T}^d \cap S_j} (g_j - h \chi_{S_j}) P_j \  d\Theta \right]^{\frac{1}{\kappa}} \notag \\
  & = C' \left[ \int_{\mathbb{T}^d} (g - h) P \  d\Theta \right]^{\frac{1}{\kappa}} \ \ (\text{as the sets } S_j \text{ are almost disjoint}),
\end{align}
whenever $\| g -h \|_{L^1(\mathbb{T}^d)} \leq \delta \bydef \min_{1\leq j\leq 4^d} C_j \frac{|q_{\boldsymbol{m}_j}|^{1/|\boldsymbol{m}_j|}}{8d |\boldsymbol{n}|}.$
The constants in the above estimation are obtained  as follows: 
$$ C= 2^d \max_{1\leq l \leq 2^d} \left[ 2^{d |\boldsymbol{n}|}C_l^{|\boldsymbol{m}_l|}(1+ |\boldsymbol{m}_l|)\right]^{\frac{1}{|\boldsymbol{m}_l|+1}}, $$
 $$\kappa =  \min_{1\leq l \leq 2^d} \frac{1}{|\boldsymbol{m}_l|+1} \leq 1, \text{ and } C'= 2^{(1-\frac{1}{\kappa})d} C. $$
Therefore, 
\begin{align} \label{Eq:combined-estimate}
\| g -h \|_{L^1(\mathbb{T}^d)}^\kappa \leq C' \left| \int_{\mathbb{T}^d} (g - h) P \  d\Theta \right|,
\end{align}
whenever $\| g -h \|_{L^1(\mathbb{T}^d)} \leq \delta$.

We would like to point out that the arguments used after the second inequality to obtain \eqref{Eq:Putinar-estimate} along with \eqref{Eq:estimate-convergence} yield
\begin{align} \label{Eq:Mainak}
            \int_{\mathbb{T}^d} (g - h) P \  d\Theta  \geq  \sum_{j=1}^{2^d} \Lambda_{|Q_j|}(\|\tilde{g}_j - v_j \|_{L^1(\lambda)}).
\end{align}
Again, 
\begin{align}\label{Eq:moment-estimate}
\left| \int_{\mathbb{T}^d} (g - h) P \  d\Theta \right| = &  \left| \sum_{\boldsymbol{\beta}\in \Gamma_{\boldsymbol{n}}} r_{\boldsymbol{\beta}}\Re \widehat{(g - h)}(\boldsymbol{\beta}) + \tilde{r}_{\boldsymbol{\beta}}\Im \widehat{(g - h)}(\boldsymbol{\beta})  \right|  \notag \\
                                                     \leq & \left[ \sum_{\boldsymbol{\beta}\in \Gamma_{\boldsymbol{n}}} (|r_{\boldsymbol{\beta}}|^2 + |\tilde{r}_{\boldsymbol{\beta}}|^2) \sum_{\boldsymbol{\beta}\in \Gamma_{\boldsymbol{n}}} 
 | \widehat{(g - h)}(\boldsymbol{\beta})|^2 \right]^{1/2}  \notag \\
                                                     \leq & C(P) \left[ \sum_{\boldsymbol{\beta}\in \Gamma_{\boldsymbol{n}}} 
 | L_{\boldsymbol{\beta}} (c_{\boldsymbol{m}}(R))_{\boldsymbol{m} \in \Gamma_{\boldsymbol{\beta}}} -  L_{\boldsymbol{\beta}} (c_{\boldsymbol{m}}(f))_{\boldsymbol{m} \in \Gamma_{\boldsymbol{\beta}}}|^2 \right]^{1/2}.
\end{align}
Recall from \eqref{Eq:univ-poly} in Subsection \ref{subsection:Exp-transform} that the functions $L_{\boldsymbol{\beta}}$ are rational having poles outside a neighbourhood $\mathcal{N}$ of $(c_{\boldsymbol{m}}(R))_{\boldsymbol{m} \in \Gamma_{\boldsymbol{n}}}$ and hence, by the mean value theorem, there exists a constant $C(R)>0$ depending on $R$ and $\mathcal{N}$ such that
\begin{align} \label{Eq:moment-conv}
\sum_{\boldsymbol{\beta}\in \Gamma_{\boldsymbol{n}}} 
 | L_{\boldsymbol{\beta}} (c_{\boldsymbol{m}}(R))_{\boldsymbol{m} \in \Gamma_{\boldsymbol{\beta}}} -  L_{\boldsymbol{\beta}} (c_{\boldsymbol{m}}(f))_{\boldsymbol{m} \in \Gamma_{\boldsymbol{\beta}}}|^2 \leq C(R) \|T_{\boldsymbol{n}}(R) - T_{\boldsymbol{n}}(f)\|_2^2,
\end{align}
whenever $(c_{\boldsymbol{m}}(f))_{\boldsymbol{m} \in \Gamma_{\boldsymbol{n}}} \in \mathcal{N}.$

\subsection*{Concluding part} Let $\{ \varepsilon_s \}_s$ be a sequence of positive real numbers and let $F>0$ $\lambda-$a.e. If $\Lambda_{F}(\varepsilon_s) \rightarrow 0$ as $s\rightarrow \infty$, then it is easy to observe from the definition of $\Lambda_{F}$ that $\varepsilon_s \rightarrow 0$ as $s\rightarrow \infty$; cf. \cite[Lemma 2.7]{Lewis}. 

Thus utilizing this observation along with \eqref{Eq:Mainak}, \eqref{Eq:moment-estimate} and \eqref{Eq:moment-conv}, we can choose an $\rho >0$ such that 
\begin{align*}
\sum_{j=1}^{2^d}  \|\tilde{g}_j - (h \chi_{S_j}) \circ \Phi_j \|_{L^1(\lambda)} \leq \frac{\delta}{2^d},
\end{align*}
whenever $ \|T_{\boldsymbol{n}}(R) - T_{\boldsymbol{n}}(f)\|_2  \leq \rho.$
Therefore, from \eqref{Eq:Mainak-2} we obtain 
\begin{align}\label{Eq:Mainak-3}\rho
\| g - h \|_{L^1(\mathbb{T}^d)} \leq \delta \ \text{  whenever } \|T_{\boldsymbol{n}}(R) - T_{\boldsymbol{n}}(f)\|_2 \leq \rho.
\end{align}

\vspace{2mm}
Now, for the Herglotz-Nevanlinna function $f$, if $\|T_{\boldsymbol{n}}(R) - T_{\boldsymbol{n}}(f)\|_2 \leq \rho$, then we take $\mathcal{N}$ to be the $\rho$-neighbourhood of $(c_{\boldsymbol{m}}(R))_{\boldsymbol{m} \in \Gamma_{\boldsymbol{n}}}$ such that \eqref{Eq:moment-conv} to hold true. Thus, combining \eqref{Eq:pointwise-estimate}, \eqref{Eq:combined-estimate}, \eqref{Eq:moment-estimate}, \eqref{Eq:moment-conv} and \eqref{Eq:Mainak-3}, we have 
\begin{align*}
 &|R(\boldsymbol{z}) - f(\boldsymbol{z})| \\
  \leq & \ \pi C(K) |\psi(0) - \varphi(0)| +   \pi C' C'(K) \left| \int_{\mathbb{T}^d} (g - h) P \  d\Theta \right|^{\frac{1}{\kappa}} \\
 \leq &  \ \pi C(K) |\log R(0) - \log f(0)| +   \pi C' C'(K) \left[\sqrt{C(R)} C(P)  \|T_{\boldsymbol{n}}(R) - T_{\boldsymbol{n}}(f)\|_2 \right]^{\frac{1}{\kappa}},
\end{align*}
for every $\boldsymbol{z} \in K$ and whenever $\|T_{\boldsymbol{n}}(R) - T_{\boldsymbol{n}}(f)\|_2 \leq \rho$.
Further, applying the mean value theorem for $\log$ in a neighbourhood of $R(0)$ in $\mathbb{C}$ (which is obtained by taking orthogonal projections of points in $\mathcal{N}$ to their first component), we get
\begin{align*}
|\log R(0) - \log f(0)| \leq C'(R)  |R(0)-f(0)| \leq C'(R) \|T_{\boldsymbol{n}}(R) - T_{\boldsymbol{n}}(f)\|_2.
\end{align*}
Finally, we have 
\begin{align*}
&|R(\boldsymbol{z}) - f(\boldsymbol{z})| \\
\leq &  \pi C(K)C'(R) \|T_{\boldsymbol{n}}(R) - T_{\boldsymbol{n}}(f)\|_2 +   \pi C' C'(K) \left[\sqrt{C(R)} C(P)  \|T_{\boldsymbol{n}}(R) - T_{\boldsymbol{n}}(f)\|_2 \right]^{\frac{1}{\kappa}}  
\end{align*}
and hence, whenever $\|T_{\boldsymbol{n}}(R) - T_{\boldsymbol{n}}(f)\|_2 \leq \rho$ (without any loss of generality we take $\rho<1$) and $\boldsymbol{z} \in K$,
\begin{align*}
&|R(\boldsymbol{z}) - f(\boldsymbol{z})|^\kappa \leq C(K, R, P) \|T_{\boldsymbol{n}}(R) - T_{\boldsymbol{n}}(f)\|_2,
\end{align*}
where $C(K, R, P) = 2^{\kappa-1} \max \{(\pi C(K)C'(R))^\kappa, ( \pi C' C'(K))^\kappa \sqrt{C(R)} C(P) \}$. In fact, the polynomial $P$ depends only on $R$ and hence we may write $C(K, R, P)$ as $C(K, R)$ and so, 
$$
\|R-f \|_{\infty, K}^\kappa \leq C(K, R) \|T_{\boldsymbol{n}}(R) - T_{\boldsymbol{n}}(f)\|_2 $$ whenever $\|T_{\boldsymbol{n}}(R) - T_{\boldsymbol{n}}(f)\|_2 \leq \rho$.
\end{proof}

\begin{remark}
	{
We have noticed in the previous section that there is a dictionary between rational inner functions and Cayley rational inner functions established by Cayley's transform. In particular, the $\boldsymbol{\alpha}$-th Taylor coefficient of a Schur function is a universal rational function depending on the Taylor coefficients up to entrywise order $\boldsymbol{\alpha}$ of its Cayley transform, and vice versa. Thus, if we consider $R$ and $f$ as appearing in Theorem \ref{Thm:supres-polydisk} with $\tilde{R}$ and $\tilde{f}$ being their inverse Cayley transforms respectively, then $\|T_{\boldsymbol{n}}(R) - T_{\boldsymbol{n}}(f) \|_2$ and $\|T_{\boldsymbol{n}}(\tilde{R}) - T_{\boldsymbol{n}}(\tilde{f}) \|_2$ are comparable with universal bounds. Moreover, note that $	\|R - f\|_{\infty, K}$ and $	\|\tilde{R} - \tilde{f}\|_{\infty, K}$ are comparable, too. Hence, a superresolution estimate for rational inner functions can be proved in the Schur class context.}
\end{remark}
\subsection{Superresolution in the Euclidean ball} Let $\mathbb{B}_d$ denote the Euclidean unit ball in $\C^d$. Not surprising, the structure of holomorphic maps $F: \mathbb{B}_d \longrightarrow \mathbb{B}_d$ is simpler and similar to the single variable case $d=1$. 
Exploiting the fact that $\mathbb{B}_d$ is a homogeneous space with a well understood group of holomorphic automorphisms,  Khudaĭberganov \cite{Khud-88} developed an analogue of Schur's algorithm for these maps. From this framework a natural superresolution theorem emerges. 

Given a holomorphic map $f: \mathbb{B}_d \rightarrow \mathbb{B}_d$ with $f=(f_1, \dots, f_d)$, we denote its Taylor series around the origin as 
$$
f(\boldsymbol{z}) = \sum_{\boldsymbol{\alpha} \in \mathbb{N}_0^d} c^{(\boldsymbol{\alpha})}(f) \boldsymbol{z}^{\boldsymbol{\alpha}} \ \ \text{ for } \boldsymbol{z} \in \mathbb{B}_d,
$$
where $c^{(\boldsymbol{\alpha})}(f) = \left( c_{\boldsymbol{\alpha}}(f_1), \dots, c_{\boldsymbol{\alpha}}(f_d)  \right)$ with $c_{\boldsymbol{\alpha}}(f_j)$ being the $\boldsymbol{\alpha}$-th Taylor coefficient of $f_j$. The affine Taylor polynomial at the origin, of $f$ is:
 $$T_{\boldsymbol{1}}(f) \bydef \sum_{|\boldsymbol{\alpha}| \leq 1} c^{(\boldsymbol{\alpha})}(f) \boldsymbol{z}^{\boldsymbol{\alpha}}.$$ Let
$$
H_j \bydef f - T_{\boldsymbol{1}}(f) = \sum_{|\boldsymbol{\alpha} | \geq 2}  c^{(\boldsymbol{\alpha})}(f) \boldsymbol{z}^{\boldsymbol{\alpha}}.
$$
The radial limits of the bounded holomorphic functions $f_j$ belong to in $L^\infty(\partial \mathbb{B}_d, \sigma)$ and we continue to denote them by the same symbol $f_j$. Here $\sigma$ is the Lebesgue surface measure on $\partial \mathbb{B}_d$. The $L^2$-norm of $f$ is given by
 \begin{align*}
  \|f\|_{L^2(\partial \mathbb{B}_d)} & \bydef \left[ \sum_{j=1}^d \|f_j\|^2_{L^2(\partial \mathbb{B}_d)} \right]^{\frac{1}{2}}
   = \left[ \sum_{j=1}^d \sum_{\boldsymbol{\alpha} \in \mathbb{N}_0^d} \| c_{(\boldsymbol{\alpha})}(f_j)\|^2 \right]^{\frac{1}{2}}
   =  \left[ \sum_{\boldsymbol{\alpha}\in \mathbb{N}_0^d} \| c^{(\boldsymbol{\alpha})}(f) \|^2 \right]^{\frac{1}{2}}.
 \end{align*}
The supperresolution result for prescribed affine data is stated as follows.
\begin{theorem}\label{Thm:Superresolution-ball}
Let $F$ be a holomorphic automorphism of $\mathbb{B}_d$ and let $f: \mathbb{B}_d \rightarrow \mathbb{B}_d$ be a holomorphic map. There exists a positive constant $C(F)$ (depending only on $F$) such that
$$\|F-f\|_{L^2(\partial \mathbb{B}_d)}^2 \leq C(F) \left[d(d+1) \rho^2  + \sqrt{d(d+1)} \rho \right] $$
where $\rho = \|T_{\boldsymbol{1}}(F) - T_{\boldsymbol{1}}(f)\|_{L^2(\partial \mathbb{B}_d)} $.
\end{theorem}
\begin{proof}
In view of the well known structure of the holomorphic automorphisms of $\mathbb{B}_d$, the map  $F$ is rational of degree one:
\begin{align*}
F(\boldsymbol{z})=(F_1(\boldsymbol{z}), \dots, F_d(\boldsymbol{z})) = \frac{a^{(0)} + a^{(1)}z_1 + \dots + a^{(d)}z_d}{b_0 + b_1 z_1 + \dots + b_d z_d}
\end{align*}
with $a^{(j)} = (a^{(j)}_1, \dots, a^{(j)}_d) \in \C^d$ and $b_j \in \C $ for $0\leq j \leq d$ and the denominator, $b(\boldsymbol{z})= b_0 + b_1 z_1 + \dots + b_d z_d$ does not vanish in $\overline{\mathbb{B}}_d$. See \cite[Chapter 2]{Rudin} for deails.

Suppose $f=(f_1, \dots, f_d)$. Then 
\begin{align}\label{Eq:Taylor-pol-dif}
F - f = T_{\boldsymbol{1}}(F-f) + \sum_{j=1}^d H_j (F-f).
\end{align}
Thus, for each $1\leq j \leq d$  and $\boldsymbol{z} \in \mathbb{B}_d$,
\begin{align*}
& (F_j - f_j)(\boldsymbol{z}) = \sum_{|\boldsymbol{\alpha}| \leq 1}  c_{\boldsymbol{\alpha}}(F_j - f_j)  \boldsymbol{z}^{\boldsymbol{\alpha}} + H_j (F-f)(\boldsymbol{z}),
\end{align*}
consequently,
\begin{align}\label{Eq:Taylor-2}
 b(\boldsymbol{z})  \left[ f_j(\boldsymbol{z})  + \sum_{|\boldsymbol{\alpha}| \leq 1}  c_{(\boldsymbol{\alpha})}(F_j - f_j)  \boldsymbol{z}^{\boldsymbol{\alpha}} \right]  = b(\boldsymbol{z})  F_j(\boldsymbol{z})  -  b(\boldsymbol{z}) H_j (F-f)  (\boldsymbol{z}). 
\end{align}

Set $$\rho = \|T_{\boldsymbol{1}}(F) - T_{\boldsymbol{1}}(f)\| = \left[\sum_{j=1}^d \sum_{|\boldsymbol{\alpha}| \leq 1} | c_{(\boldsymbol{\alpha})}(F_j - f_j)|^2 \right]^{1/2}.$$
From $$\sum_{j=1}^d |f_j(\boldsymbol{\xi})|^2 \leq 1, \ \   \sigma-\text{a.e. on } \partial \mathbb{B}_d $$ and 
 $$\sum_{j=1}^d  \sum_{|\boldsymbol{\alpha}| \leq 1} |\boldsymbol{\xi}^{\boldsymbol{\alpha}} |^2 = 2 \ \ \text{ on } \partial \mathbb{B}_d,$$ we infer that
\begin{align}\label{Eq:Taylor-3}
&\sum_{j=1}^d \int_{\partial \mathbb{B}_d } | b(\boldsymbol{\xi})|^2  | f_j(\boldsymbol{\xi})  + \sum_{|\boldsymbol{\alpha}| \leq 1}  c_{\boldsymbol{\alpha}}(F_j - f_j)  \boldsymbol{\xi}^{\boldsymbol{\alpha}} |^2 d\sigma(\boldsymbol{\xi}) \notag \\
 &\leq  \sum_{j=1}^d \int_{\partial \mathbb{B}_d } | b(\boldsymbol{\xi})|^2 \left( | f_j(\boldsymbol{\xi})|^2 + \left| \sum_{|\boldsymbol{\alpha}| \leq 1}  c_{\boldsymbol{\alpha}}(F_j - f_j)  \boldsymbol{\xi}^{\boldsymbol{\alpha}} \right|^2 + 2 \left| \sum_{|\boldsymbol{\alpha}| \leq 1}  c_{\boldsymbol{\alpha}}(F_j - f_j)  \boldsymbol{\xi}^{\boldsymbol{\alpha}} \right| \right)d\sigma(\boldsymbol{\xi}) \notag \\
 & \leq \int_{\partial \mathbb{B}_d } | b(\boldsymbol{\xi})|^2 \left( 1+ \rho^2 \sum_{j=1}^d  \sum_{|\boldsymbol{\alpha}| \leq 1} |\boldsymbol{\xi}^{\boldsymbol{\alpha}} |^2 + 2 \rho  \left[ \sum_{j=1}^d  \sum_{|\boldsymbol{\alpha}| \leq 1} |\boldsymbol{\xi}^{\boldsymbol{\alpha}} |^2 \right]^{\frac{1}{2}} \right)d\sigma(\boldsymbol{\xi}) \notag \\
 & \leq  \left( 1+ d(d+1) \rho^2  + \sqrt{d(d+1)} \rho  \right) \|b\|_{L^2(\partial \mathbb{B}_d)}^2.
\end{align}
 
Again, $F$ being an automorphism we have $F(\partial \mathbb{B}_d ) = \partial \mathbb{B}_d $ and therefore,
\begin{align}\label{Eq:Taylor-4}
\sum_{j=1}^d |b(\boldsymbol{\xi}) F_j(\boldsymbol{\xi}) |^2 = |b(\boldsymbol{\xi})|^2 \ \ \text{ on }\partial \mathbb{B}_d.
\end{align}
Now \eqref{Eq:Taylor-2} and \eqref{Eq:Taylor-3} yield:
\begin{align}\label{Eq:Taylor-5}
\sum_{j=1}^d \|  b F_j  -  b H_j (F-f) \|_{L^2(\partial \mathbb{D}^d)}^2 \leq \left( 1+ d(d+1) \rho^2  + \sqrt{d(d+1)} \rho  \right) \|b\|_{L^2(\partial \mathbb{B}_d)}^2.
\end{align} 
Note that, $b(\boldsymbol{z}) F_j(\boldsymbol{z}) = a^{(j)}_0 + a^{(j)}_1z_1 + \dots + a^{(j)}_d z_d $ is a linear polynomial and each term of $b(\boldsymbol{z}) H_j (F-f)(\boldsymbol{z})$ has total degree greater than or equal to $2$. Thus $bF_j$ and $b H_j (F-f)$ are orthogonal in $L^2(\partial \mathbb{B}_d)$. Hence, by \eqref{Eq:Taylor-4} and \eqref{Eq:Taylor-5}, we get
\begin{align*}
&\sum_{j=1}^d \|  b F_j \|_{L^2(\partial \mathbb{D}^d)}^2  +  \| b H_j (F-f) \|_{L^2(\partial \mathbb{D}^d)}^2 \\
&= \|b\|_{L^2(\partial \mathbb{B}_d)}^2 + \sum_{j=1}^d  \| b H_j (F-f) \|_{L^2(\partial \mathbb{D}^d)}^2 \\
& \leq   \left( 1+ d(d+1) \rho^2  + \sqrt{d(d+1)} \rho  \right) \|b\|_{L^2(\partial \mathbb{B}_d)}^2.
\end{align*}
Therefore,
\begin{align*}
b_{\text{min}} \sum_{j=1}^d \| H_j (F-f) \|_{L^2(\partial \mathbb{D}^d)}^2 & \leq  \sum_{j=1}^d  \| b H_j (F-f) \|_{L^2(\partial \mathbb{D}^d)}^2 \\
                                                                           & \leq (d(d+1) \rho^2  + \sqrt{d(d+1)} \rho ) \|b\|_{L^2(\partial \mathbb{B}_d)}^2,
\end{align*} 
where the minimum of $|b|$ over $\overline{\mathbb{B}}_d$, $b_{\text{min}}$ is positive as $b$ is non-vanishing on $\overline{\mathbb{B}}_d$. Finally, we obtain from \eqref{Eq:Taylor-pol-dif} 
$$
\|F-f\|_{L^2(\partial \mathbb{B}_d)}^2  \leq \frac{1}{b_{\text{min}}} \left[d(d+1) \rho^2  + \sqrt{d(d+1)} \rho \right] \|b\|_{L^2(\partial \mathbb{B}_d)}^2,
                                 $$
with the constant $C(F)= \frac{ \|b\|_{L^2(\partial \mathbb{B}_d)}^2}{b_{\text{min}}}$ which only depends on the denominator $b$ of $F$. This completes the proof.
\end{proof}

Regardless to say that any similar superresolution statement referring to a higher degree interpolatory data, containing as a subset the affine part of the Taylor polynomial at the origin, can be derived 
from Theorem \ref{Thm:Superresolution-ball}. To be more precise, the distance between jets of degree one of two functions is always less than or equal than the distance between their jets of higher degree.

\subsection{Final comments} \label{Exp-open-questions}

So far, we have exploited the fact that rational inner functions in the polydisk are determined by a section of their Taylor expansion at zero, defined by a top leading monomial multi-index.
We can relax the set of coefficients constraints and seek corresponding extremal elementary functions which are determined by a finite section of their Taylor expansion.

The simplest case of functions $f : \D^d \longrightarrow \D$ with a prescribed degree one Taylor expansion was analyzed by Dautov and Khudaĭberganov \cite{Dautov-Khud}.
In that case, the uniqueness of the rational extremal function, essential for our superresolution theorem, fails. Examples are elementary and illustrate the lack of a Schwarz Lemma 
in the multivariable setting. A far reaching study of this sharp departure from the single variable situation goes back to H. Cartan \cite{Cartan}.

To be more specific, we consider the families of functions:
$$ f_\lambda(z_1, z_2) = z_1 + \lambda z_2^2 \ \ \text{ for }(z_1, z_2)\in \mathbb{B}_2 $$ 
and 
$$  
g_\lambda(z_1,z_2) = \frac{z_1+ z_2}{2} + \lambda (z_1^2- z_2^2) \ \ \text{ for }(z_1, z_2)\in \D^2. 
$$ 
They have prescribed affine Taylor polynomials $z_1$, respectively $\frac{z_1 + z_2}{2}$ and for $\lambda \in [0, 1/4]$
they satisfy $\| f_\lambda \|_{\infty, \mathbb{B}_2} = 1, \ \| g_\lambda \|_{\infty, \D^2} = 1.$ The family $\{f_\lambda\}$ appeared in \cite{Dautov-Khud}. A Cayley transform carries over the second counter-example to the Herglotz-Nevanlinna class on $\D^2$.

In spite of relentless attempts aimed at understanding and classifying unique solutions to the finite point Hermite bounded holomorphic interpolation in classical domains, only very particular low degree rational functions were identified, see for instance \cite{Agler-McCarthy, 
Kosinski-Zwonek}. By Hermite interpolation we mean finitely many points with prescribed finite order jets (which means Taylor polynomials of prescribed degrees around those points). It is clear from the cited works that invariant metrics and the explicit description of holomorphic automorphisms of the respective domains enter into discussion. One step further, we may ask the natural question: {\it is there a superresolution estimate in a neighbourhood of each of these special uniquely determined finite interpolation data solutions?}


\begin{thebibliography}{10}

\bibitem{Agler-McCarthy} Agler, Jim; McCarthy, John E. {\em The three point Pick problem on the bidisk.} New York J. Math. 6 (2000), 227--236.

\bibitem{Agler-McCarthy-Young}
Agler, Jim; McCarthy, John Edward; Young, Nicholas. {\em Operator analysis—Hilbert space methods in complex analysis}. Cambridge Tracts in Mathematics, 219. Cambridge University Press, Cambridge, 2020. xv+375 pp. 

\bibitem{Ahiezer-Krein}
Aheizer, N. I.; Krein, M.
\newblock {\em Some questions in the theory of moments}.
\newblock translated by W. Fleming and D. Prill. Translations of Mathematical
  Monographs, Vol. 2. American Mathematical Society, Providence, R.I., 1962.
  
  \bibitem{Aleksandrov}
  Aleksandrov, A. B. {\em Function theory in the ball.} (Russian) Current problems in mathematics. Fundamental directions, Vol. 8, 115--190, 274, Itogi Nauki i Tekhniki, Akad. Nauk SSSR, Vsesoyuz. Inst. Nauchn. i Tekhn. Inform., Moscow, 1985.
  
  \bibitem{Bhowmik-BLMS} Bhattacharyya, Tirthankar; Bhowmik, Mainak; Kumar, Poornendu. {\em Herglotz's representation and Carathéodory's approximation}. Bull. Lond. Math. Soc. 56 (2024), no. 12, 3752–3776.
  
  \bibitem{Bhowmik-Kumar-JFA} Bhowmik, Mainak; Kumar, Poornendu. {\em Function theory on quotient domains related to the polydisc}.
J. Funct. Anal. 289 (2025), no. 6, Paper No. 110978, 29 pp.

\bibitem{BK-PAMS2}
Bhowmik, Mainak; Kumar, Poornendu.
 {\em Herglotz representation for operator-valued function on a set associated with test functions}.
Proc. Amer. Math. Soc. 153 (2025), no. 7, 2949–2963.
  
 \bibitem{Budisic-Putinar}
  Budišić, Marko; Putinar, Mihai. {\em Conditioning moments of singular measures for entropy optimization. I}. Indag. Math. (N.S.) 23 (2012), no. 4, 848--883.


\bibitem{Candes-Fernandez}
Candès, Emmanuel J.; Fernandez-Granda, Carlos. {\em Towards a mathematical theory of super-resolution}. Comm. Pure Appl. Math. 67 (2014), no. 6, 906--956. 


\bibitem{Cartan}
Cartan, H. {\em Sur les fonctions de deux variables complexes. Les transformations d'un domaine borné $D$ en un domaine intérieur à $D$}. (French) Bull. Soc. Math. France 58 (1930), 199--219.


\bibitem{Dautov-Khud}
Dautov, Sh. A.; Khudaĭberganov, G. {\em The Carathéodory-Fejér problem in higher-dimensional complex analysis}. (Russian) Sibirsk. Mat. Zh. 23 (1982), no. 2, 58--64, 215.



\bibitem{Dieu-2018}
Dieu, Thieu~Anh.
\newblock {\em Volume estimates of sublevel sets of real polynomials}.
\newblock { Ann. Polon. Math.}, 121(2):157--174, 2018.


\bibitem{Donoho}
Donoho, David L. {\em Superresolution via sparsity constraints}. SIAM J. Math. Anal. 23 (1992), no. 5, 1309--1331.

\bibitem{Eschmeier-Patton-Putinar} Eschmeier, Jörg; Patton, Linda; Putinar, Mihai. {Carathéodory-Fejér interpolation on polydisks}. Math. Res. Lett. 7 (2000), no. 1, 25–34.

\bibitem{FF} Foias, Ciprian; Frazho, Arthur E. {\em The commutant lifting approach to interpolation problems}. Operator Theory: Advances and Applications, 44. Birkhäuser Verlag, Basel, 1990. {\rm xxiv}+632 pp. ISBN: 3-7643-2461-9

\bibitem{Gamboa}
Gamboa, F; Gassiat, E.
\newblock {\em Sets of superresolution and the maximum entropy method on the mean}.
\newblock { SIAM J. Math. Anal.}, 27(4):1129--1152, 1996.

\bibitem{Khrushchev}
Khrushchev, Sergey. {\em Orthogonal polynomials and continued fractions. From Euler's point of view}. Encyclopedia of Mathematics and its Applications, 122. Cambridge University Press, Cambridge, 2008. xvi+478 pp.

\bibitem{Khud-88}
Khudaĭberganov, G. {\em The Carathéodory-Fejér problem in the generalized unit disk}. (Russian) ; translated from Sibirsk. Mat. Zh. 29 (1988), no. 6, 160--166 Siberian Math. J. 29 (1988), no. 6, 1000--1005 (1989) 



\bibitem{Knese}
Knese, G.  {\em Rational Inner Functions on the Polydisk: A Survey}. In: Alpay, D., Sabadini, I., Colombo, F. (eds) Operator Theory. Springer, Basel. 25 pp., 2025.

\bibitem{KP}
Kor\'anyi, A.; Puk\'anszky, L. {\em Holomorphic functions with positive real part on polycylinders}. Trans. Amer. Math. Soc. 108 (1963), 449--456.

\bibitem{Koranyi-Vagi} Kor\'anyi, Adam; Vági, Stephen. {\em Rational inner functions on bounded symmetric domains}. Trans. Amer. Math. Soc. 254 (1979), 179–193.

\bibitem{Kosinski-Zwonek}
Kosiński, Łukasz; Zwonek, Włodzimierz. {\em Nevanlinna-Pick problem and uniqueness of left inverses in convex domains, symmetrized bidisc and tetrablock.} J. Geom. Anal. 26 (2016), no. 3, 1863--1890.


\bibitem{Krein-Nudelman}
Kre\u{\i}n, M.~G.; Nudelman, A.~A.
\newblock {\em The {M}arkov moment problem and extremal problems}.
\newblock American Mathematical Society, Providence, R.I., 1977.
\newblock Ideas and problems of P. L. \v{C}eby\v{s}ev and A. A. Markov and
  their further development, Translated from the Russian by D. Louvish,
  Translations of Mathematical Monographs, Vol. 50.


\bibitem{Lewis}
Lewis, A.~S. 
\newblock {\em Superresolution in the {M}arkov moment problem}.
\newblock { J. Math. Anal. Appl.}, 197(3):774--780, 1996.


\bibitem{Mallat}
Mallat, Stéphane; Yu, Guoshen.
{\em Super-resolution with sparse mixing estimators}. 
IEEE Trans. Image Process. 19, No. 11, 2889-2900 (2010).


\bibitem{Nguyen-Milanfar-Golub}
Nguyen, Nhat; Milanfar, Peyman; Golub, Gene.
{\em A computationally efficient superresolution image reconstruction algorithm}.
IEEE Trans. Image Process. 10, No. 4, 573-583 (2001).


\bibitem{Pfister}
Pfister, Albrecht. {\em Über das Koeffizientenproblem der beschränkten Funktionen von zwei Veränderlichen}. (German) Math. Ann. 146 (1962), 249--262.

\bibitem{Putinar-JAT}
Putinar, Mihai. {\em Extremal solutions of the two-dimensional  L -problem of moments. II}.
J. Approx. Theory 92 (1998), no. 1, 38–58.

\bibitem{Putinar}
Putinar, Mihai. {\em Superresolution of principal semi-algebraic sets}. Anal. Math. Phys. 11 (2021), no. 3, Paper No. 107, 14 pp.

\bibitem{Rudin}
 Rudin, Walter.
 {\em Function Theory in the Unit Ball of $\mathbb{C}^n$}. Reprint of the 1980 edition, Classics in Mathematics, Springer-Verlag, Berlin, ISBN 978-3-540-68272-1, 2008, xiv+436 pp.

\bibitem{Schur} Schur, Isaac. {\em \"Uber Potenzreihen, die im Innern des Einheitskreises beschr\"ankt sind}.
J. Math. 147, 205-232 (1917), 148, 122-145 (1918).

\bibitem{Varchenko-1976}
 Var\v{c}enko, A.~N.
\newblock {\em Newton polyhedra and estimates of oscillatory integrals}.
\newblock { Funkcional. Anal. i Prilo\v{z}en.}, 10(3):13--38, 1976.


\bibitem{Woerdeman} Woerdeman, Hugo J. {\em Positive Carathéodory interpolation on the polydisc}. Integral Equations Operator Theory 42 (2002), no. 2, 229--242.


\end{thebibliography}
\end{document}